\documentclass[11pt]{amsart}
\usepackage[T1]{fontenc}
\usepackage[utf8]{inputenc}
\usepackage{lmodern,microtype}
\usepackage{amsmath,amssymb,mathtools}
\usepackage{enumitem,needspace,booktabs,array,tikz-cd}
\usepackage{xurl,xcolor,float}
\usepackage[letterpaper,margin=1.02in]{geometry}
\usepackage{hyperref}
\hypersetup{colorlinks=true,linkcolor=blue!45!black,citecolor=blue!45!black,urlcolor=blue!45!black,pdftitle={Thom morphisms, self-intersection, and the sectional category of a section},pdfauthor={Paul-Eugène Parent}}
\allowdisplaybreaks[2]
\usepackage{etoolbox}
\makeatletter
\patchcmd{\@maketitle}{\global\topskip42\p@\relax}{\global\topskip\normaltopskip\relax}{}{}
\makeatother
\apptocmd{\thebibliography}{\setlength{\itemsep}{0pt plus 4pt}}{}{}

\newcommand{\Q}{\mathbb Q}
\newcommand{\TC}{\operatorname{TC}}
\newcommand{\MTC}{\operatorname{MTC}}
\newcommand{\secat}{\operatorname{secat}}
\newcommand{\relcat}{\operatorname{relcat}}
\newcommand{\msecat}{\operatorname{Msecat}}
\newcommand{\id}{\mathrm{id}}
\newtheorem{theorem}{Theorem}[section]
\newtheorem{proposition}[theorem]{Proposition}
\newtheorem{lemma}[theorem]{Lemma}
\newtheorem{corollary}[theorem]{Corollary}
\newtheorem{maintheorem}{Theorem}

\theoremstyle{definition}
\newtheorem{definition}[theorem]{Definition}
\newtheorem{example}[theorem]{Example}
\theoremstyle{remark}
\newtheorem{remark}[theorem]{Remark}
\title[Thom morphisms and sectional category]{Thom morphisms, self-intersection, and the sectional category of a section}
\DeclareRobustCommand{\fundingmark}{\textsuperscript{*}}
\author[Paul-Eugène Parent]{Paul-Eugène Parent\fundingmark}
\thanks{\makebox[0pt][r]{\fundingmark\enspace}This research was supported by the Natural Sciences and Engineering Research Council of Canada (NSERC), grant no.~149521.}
\address{Department of Mathematics and Statistics, Faculty of Science,
University of Ottawa, STEM Complex, Room 361, 150 Louis-Pasteur Private,
Ottawa, Ontario K1N 6N5, Canada}
\email{pparent@uottawa.ca}
\date{November 9th, 2026}
\subjclass[2020]{Primary 55M30; Secondary 55P62, 13D02}
\keywords{Sectional category, relative category, Thom morphism, self-intersection, Poincar\'e duality, Sullivan model}
\begin{document}
\begin{abstract}
For sections between simply connected spaces of finite rational type,
we compare module sectional category, sectional category and relative
category after rationalization. Fibre Poincar\'e duality turns the
compatibility of derived evaluation with the section into a derived
self-intersection formula. We prove equality of the three invariants
under any of four hypotheses: a base-dimension bound, finite total rational
homotopy of the fibre, a smooth section with closed manifold fibre,
or a compatible relative Frobenius model. Projective bundles over
$\mathbb{CP}^k$ with fibre $\mathbb{CP}^r$ have common value $k+r$,
which is unbounded even for fixed fibre. Whether fibre Poincar\'e
duality alone suffices remains open.
\end{abstract}
\maketitle
\section{Introduction}
\label{sec:introduction}

A section $f:X\to Y$ of a fibration $g:Y\to X$ distinguishes a copy
of the base inside the total space. How many local homotopy sections
of this inclusion are needed to cover $Y$? Sectional category answers
this question; relative category asks those sections to respect the
distinguished copy of $X$. Their rational versions have an algebraic
lower bound, module sectional category, obtained by asking for
retractions of modules instead of algebras:
\[
 \msecat(f_\Q)\leq\secat(f_\Q)\leq\relcat(f_\Q).
\]
We ask when the geometry of the fibre of $g$ forces equality.

The classical self-intersection formula supplies the guiding idea.
For a smooth section with oriented normal bundle $\nu$, the Thom
map satisfies $f^*f_!(1)=e(\nu)$. We prove a derived version of this
formula for rationally retractive maps: a map $f$ is rationally
retractive when $f_\Q$ has a homotopy retraction. All spaces below
are simply connected and of finite rational type, meaning finite
dimensional rational cohomology in each degree.

\begin{maintheorem}[The Thom criterion]\label{thm:intro-thom}
Let $f:X\to Y$ be rationally retractive, and suppose the fibre of a
chosen retraction has Poincar\'e duality cohomology of formal dimension
$N$. In a split Sullivan model $\varphi:R\to A$ of the section, put
$I=\ker\varphi$. Choose a Thom morphism $\tau:A[-N]\to R$, and let
$e\in H^N(A)$ be its restriction to the section. The following are
equivalent:
\begin{enumerate}[label=\textup{(\roman*)},leftmargin=*,itemsep=2pt]
\item Canonical derived evaluation
$I\otimes_R^{\mathbf L}\mathbf R\operatorname{Hom}_R(I,R)\to R$
factors through $I$.
\item Thom multiplication
$1_I\otimes\tau:I\otimes_R^{\mathbf L}A[-N]\to I$ is zero.
\item The derived self-intersection formula holds:
\[
 \varphi\tau=\varphi^*(m_e)\qquad\text{in }\mathcal D(R).
\]
\end{enumerate}
Each condition gives
$\msecat(f_\Q)=\secat(f_\Q)=\relcat(f_\Q)$, including $\infty$.
\end{maintheorem}

Here $A$ models the base $X$, $R$ models $Y$, $m_e$ is multiplication
by $e$, and $\varphi^*$ restricts scalars to the derived category
$\mathcal D(R)$. The definitions are given in
Sections~\ref{sec:preliminaries}--\ref{sec:pd-fibres}.
With the geometric normalization of the Thom morphism, $e=e(\nu)$. In general, restricting to the base
always defines $e$; the theorem asks for an identity before that
restriction. Fibre duality supplies biduality, while self-intersection
is the extra compatibility needed for evaluation to take values in $I$.

\begin{maintheorem}[Four sufficient conditions]\label{thm:intro-equality}
Let $f:X\to Y$ be rationally retractive, with chosen retraction of
$f_\Q$ and homotopy fibre $F$. Then
$\msecat(f_\Q)=\secat(f_\Q)=\relcat(f_\Q)$ if any one of the following holds.
\begin{enumerate}[label=\textup{(\arabic*)},leftmargin=*,itemsep=3pt]
\item $H^*(F;\Q)$ is a Poincar\'e duality algebra of formal dimension
$N$, the relative generators in the split minimal model have degrees
at least $c$, and $H^{>n}(X;\Q)=0$ for integers $c\geq2$ and
$n\geq0$ satisfying $n<N+c-1$.
\item $F$ has finite total rational homotopy:
$\dim_\Q(\pi_*(F)\otimes\Q)<\infty$.
\item $f$ is a smooth section of a smooth bundle over a simply
connected smooth manifold without boundary, with closed simply
connected manifold fibre.
\item The split Sullivan model fits into a strictly commutative
square whose columns are quasi-isomorphisms and whose bottom row
is the augmentation of a relative Frobenius model
(Definition~\ref{def:relative-frobenius} and
Corollary~\ref{cor:relative-frobenius-equality}).
\end{enumerate}
\end{maintheorem}

The four proofs use, respectively, degrees, a Thom cocycle, a tubular
neighbourhood, and a perfect relative multiplication pairing. The
last route reaches the evaluation criterion directly, without a Thom
calculation. Whether it proves equality beyond the other available
tests remains open; the examples here illustrate overlapping hypotheses.

For the higher diagonal $\Delta_r:X\to X^r$, $r\geq2$, sectional
category is higher topological complexity:
$\TC_r(X)=\secat(\Delta_r)$. If $H^*(X;\Q)$ has Poincar\'e duality
of formal dimension $n$, the first projection has fibre $X^{r-1}$
of dimension $(r-1)n$, so condition~\textup{(1)} applies. Thus
\[
 \MTC_r(X)=\TC_r(X_\Q)=\relcat((\Delta_r)_\Q),
 \qquad \MTC_r(X)=\msecat((\Delta_r)_\Q).
\]
This recovers the higher-diagonal equality of
\cite[Corollary~5.2]{ParentI} by fibre integration and a degree bound.

The projective family shows both the answer and its scale. If
$L\to\mathbb{CP}^k$ is the hyperplane line bundle, a trivial line
defines a section $f$ of $\mathbb P(\underline{\mathbb C}^{\,r}\oplus L)$,
and we prove
\[
 \msecat(f_\Q)=\secat(f_\Q)=\relcat(f_\Q)=k+r\qquad(k,r\geq1).
\]
The value is unbounded with fixed fibre $\mathbb{CP}^r$. In cohomology,
the class $x$ vanishing on the section satisfies $x^{r+1}+ux^r=0$,
where $u=c_1(L)$; this relation records the twisting and explains the
growth. A second family has value $k+3$ with one fixed rationally
hyperbolic manifold fibre.

\begin{figure}[H]
\centering
\begin{tikzpicture}[x=1pt,y=1pt,>=stealth,
 roadmap/.style={draw=black!45,rounded corners=2pt,fill=black!2,
 align=center,text width=100pt,inner sep=5pt,font=\footnotesize}]
\node[roadmap,text width=260pt] (evaluation) at (0,0)
 {Biduality and ideal-valued evaluation\\
 Theorem~\ref{thm:imported-evaluation} (from \cite{ParentI})};
\node[roadmap,text width=200pt] (thom) at (-58,-70)
 {Thom criterion (Theorem~\ref{thm:thom-evaluation})\\
 Self-intersection detects ideal-valued evaluation};
\node[roadmap] (degree) at (-174,-170)
 {Base-dimension bound\\Theorem~\ref{thm:pd-fibre-dimension}\\
 Degree vanishing\\Poincar\'e diagonals};
\node[roadmap] (finite) at (-58,-170)
 {Finite total rational homotopy of the fibre\\
 Theorem~\ref{thm:finite-homotopy-fibre}\\Thom cocycle\\Homogeneous fibres};
\node[roadmap] (smooth) at (58,-170)
 {Smooth section\\Theorem~\ref{thm:smooth-section}\\
 Tubular neighbourhood\\Fixed hyperbolic fibre};
\node[roadmap] (frobenius) at (174,-170)
 {Relative Frobenius model\\Corollary~\ref{cor:relative-frobenius-equality}\\
 Double annihilator\\Projective bundles};
\draw[->] (evaluation.south) -- (thom.north);
\draw[->] (thom.south) -- (degree.north);
\draw[->] (thom.south) -- (finite.north);
\draw[->] (thom.south) -- (smooth.north);
\draw[->] (evaluation.east) -- (frobenius.north);
\end{tikzpicture}
\caption{The main proof dependencies. Each lower box gives
$\msecat(f_\Q)=\secat(f_\Q)=\relcat(f_\Q)$ under its stated
hypotheses. The Frobenius route bypasses the Thom criterion.}
\label{fig:proof-roadmap}
\end{figure}
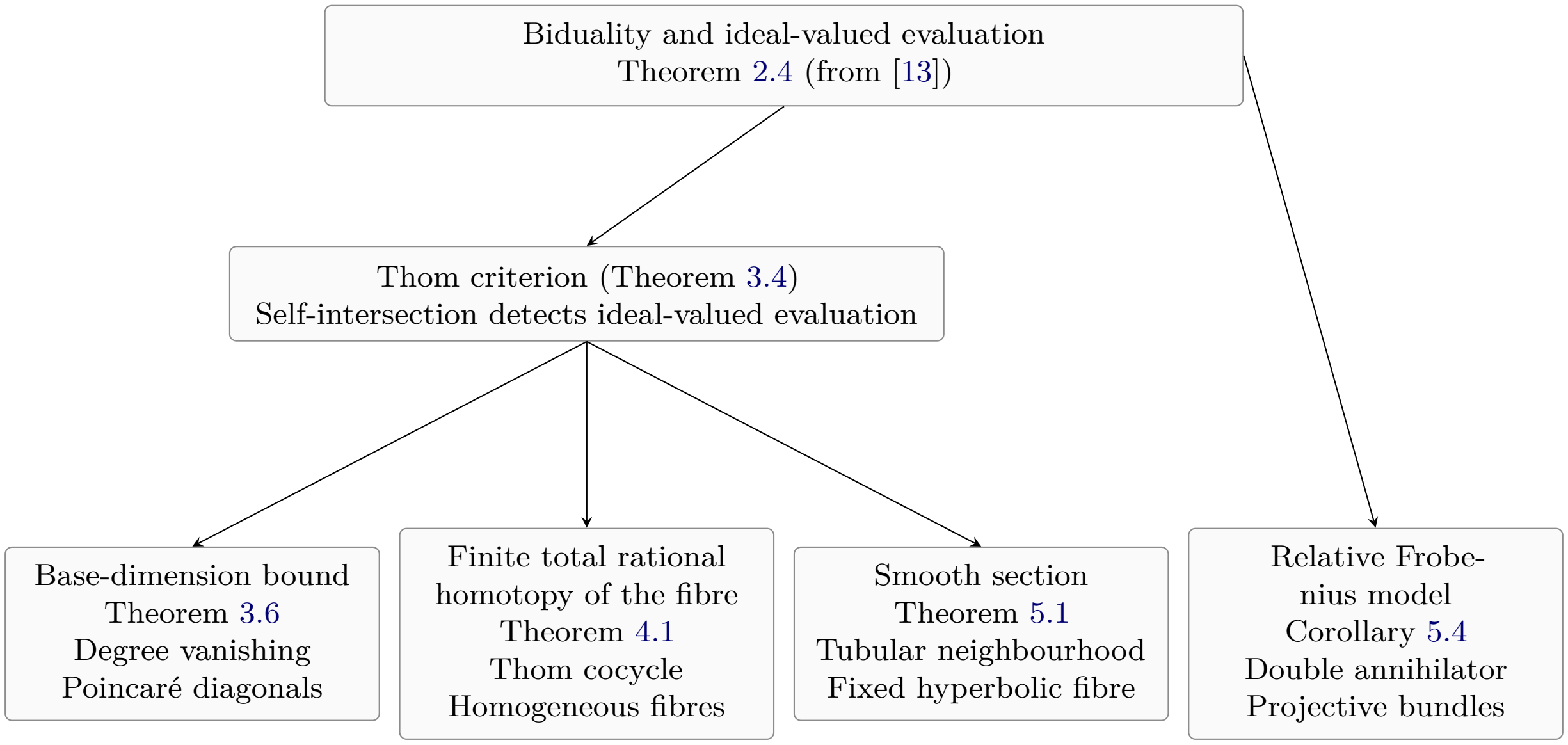

Figure~\ref{fig:proof-roadmap} records the four routes. After the
preliminaries, Section~\ref{sec:pd-fibres} proves Theorem~A and the
base-dimension case of Theorem~B. Section~\ref{sec:finite-homotopy}
constructs the Thom cocycle and treats homogeneous fibres.
Section~\ref{sec:smooth-frobenius} gives the smooth and Frobenius
arguments, and Section~\ref{sec:geometric-examples} computes the two
families. Section~\ref{sec:limits} interprets a known strict example
$\msecat<\secat$ from \cite{Parent2026} and states the open questions.
The appendices contain the relative-duality proof and the filtered
comparison for the Thom cocycle.

\section{Preliminaries}\label{sec:preliminaries}

\subsection{Conventions and sectional invariants}

The section and its retraction play different roles in the argument:
we measure the section and impose hypotheses on the fibre of the
retraction. The following conventions keep those roles fixed.
A CDGA is a commutative differential graded algebra over $\Q$.

All spaces have the homotopy type of simply connected CW complexes and
are of \emph{finite rational type}: $H^j(X;\Q)$ is finite dimensional
for every $j$. In this simply connected setting, this is equivalent to
$\pi_j(X)\otimes\Q$ being finite dimensional for every $j$.
No finiteness of the total rational homotopy is assumed.
All invariants are normalized, take values in
$\mathbb N\cup\{\infty\}$, and are applied to the rationalized section.
A map $f:X\to Y$ is \emph{rationally retractive} if there is a map
$g:Y_\Q\to X_\Q$ with $gf_\Q\simeq\id$. Replacing $g$ by a fibration
makes $f_\Q$ its homotopy section. Throughout, $F$ is the homotopy
fibre of this retraction $g$, not the fibre of $f$.

Relative Sullivan factorizations and their realization comparison
apply degree by degree under this convention
\cite[Sections~13\textup{(b)}, 14\textup{(a)}, and Proposition~15.8]{FHT2001}.

We use cohomological differentials, $M[k]^i=M^{i+k}$, and Koszul
signs. The suspension $sM$ satisfies $|sm|=|m|-1$ and
$d(sm)=-s(dm)$. The free graded-commutative algebra $\Lambda W$ is polynomial
on even generators and exterior on odd generators. A \emph{semifree}
module is free as a graded module, with an ordered filtration of its
basis (its \emph{cells}) such that each cell differential uses earlier
cells. For a homogeneous $R$-linear map $v$ of degree $k$, our conventions
are $v(rm)=(-1)^{k|r|}rv(m)$ and
$\partial v=dv-(-1)^kvd$.
Semifree resolutions compute derived tensor products and Hom.
We put $\mathbb D_RM=\mathbf R\operatorname{Hom}_R(M,R)$.
A module is \emph{perfect} if it is a retract in the derived category
of a finite semifree module.

\begin{definition}[The three sectional invariants]\label{def:three-invariants}
The \emph{sectional category} $\secat(f)$ is the least $n$ such that
$Y$ has an open cover by $n+1$ sets on which $f$ admits homotopy
sections. Equivalently, the $(n+1)$-fold fibrewise join of a fibration
replacement of $f$, denoted $g_n:G_n(f)\to Y$, admits a homotopy
section. There is a canonical map $\alpha_n:X\to G_n(f)$ with
$g_n\alpha_n\simeq f$. The \emph{relative category} $\relcat(f)$ is
the least $n$ for which a homotopy section $v$ of $g_n$ can be chosen
with $vf\simeq\alpha_n$.

Let $R$ be a Sullivan model of $Y_\Q$, and let
$R\hookrightarrow R\otimes\Lambda Z_n$ model $g_n$ for $f_\Q$.
The \emph{module sectional category} $\msecat(f_\Q)$ is the least
$n$ for which this inclusion has a differential $R$-module retraction,
or equivalently a left inverse in $\mathcal D(R)$, since the model
is semifree as an $R$-module. This condition is independent of the
chosen relative Sullivan model. Each invariant is
infinite if no such $n$ exists.
\end{definition}
These definitions and their model independence are discussed in
\cite{DEH2013,FGKV2006}; forgetting structure gives
\begin{equation}\label{eq:three-inequalities}
 \msecat(f_\Q)\leq\secat(f_\Q)\leq\relcat(f_\Q).
\end{equation}

\subsection{The split Sullivan model}

A rational retraction admits a split minimal model
\begin{equation}\label{eq:split-model}
 A\xrightarrow{\iota_A}R=A\otimes\Lambda W
       \xrightarrow{\varphi}A,
 \qquad \varphi\iota_A=\id,
 \quad \varphi(W)=0,
 \quad I=(W).
\end{equation}
Here $A$ models $X$, $R$ models $Y$, $\iota_A$ models $g$, and
$(\Lambda W,\overline d)$ models $F$. The generators have degrees at
least two and finite dimension in each degree. Minimality and the
splitting give $dW\subset I$ and
$I^m=A\otimes\Lambda^{\geq m}W$. We call word length in $W$ the
\emph{relative length}; coefficients from $A$ do not contribute to it.
One obtains this model by choosing minimal representatives of $f$ and
$g$, correcting their composite to the identity, and replacing each
relative generator $w$ by $w-\iota_A\varphi(w)$. This triangular change
preserves minimality and makes the kernel precisely $(W)$; see
\cite[Lemma~2.7]{ParentI}.

\subsection{The obstruction and evaluation criteria}

The first input places the two invariants in the same derived category $\mathcal D(R)$. It
combines the module quotient criterion of Carrasquel-Vera
\cite{CV2017} with the ideal-valued criterion of \cite[Theorem~3.1]{ParentI}.

\begin{theorem}[The obstruction criterion]
\label{thm:imported-obstruction}
Let $f:X\to Y$ be rationally retractive between simply connected spaces
of finite rational type, and use the split minimal model
\eqref{eq:split-model}. For $m\geq1$ define
\[
 o_m(f)=[I^m\hookrightarrow I]
       \in\operatorname{Hom}_{\mathcal D(R)}(I^m,I),
 \qquad \iota:I\hookrightarrow R.
\]
For every $n\geq0$,
\[
 \msecat(f_\Q)\leq n\iff\iota_*o_{n+1}(f)=0,
 \qquad
 \relcat(f_\Q)\leq n\iff o_{n+1}(f)=0.
\]
Also $\secat(f_\Q)\leq n$ if and only if
$R\to R/I^{n+1}$ admits a homotopy retraction of CDGAs.
\end{theorem}

We will also use the following elementary form of the quotient
criterion, which replaces a homotopy retraction by a strict one.

\begin{lemma}[Strict retractions]\label{lem:strict-retraction-representative}
Let $R$ be a Sullivan algebra and $j:R\hookrightarrow R'$ a
relative Sullivan inclusion. If $j$ has a left inverse in the CDGA
homotopy category, it has a strict CDGA retraction.
\end{lemma}
\begin{proof}
Since $R'$ is cofibrant and $R$ is fibrant,
represent the inverse by $a:R'\to R$. Choose a Sullivan homotopy
$H:R\to R\otimes\Lambda(t,dt)$ from $aj$ to $\id_R$.
The square with top map $H$, bottom map $a$, left map $j$, and right
map $\operatorname{ev}_0:R\otimes\Lambda(t,dt)\to R$ commutes.
The right map is a surjective quasi-isomorphism, so the relative
Sullivan lifting property gives $\widetilde H:R'\to R\otimes\Lambda(t,dt)$
with $\widetilde H j=H$ and $\operatorname{ev}_0\widetilde H=a$.
Then $\operatorname{ev}_1\widetilde H$ is the required strict retraction.
\end{proof}

A nullhomotopy on a semifree resolution $P\to I^m$ represents the
first vanishing when its values lie in $R$, and the second when its
values lie in $I$. The difference is the possible component $\varphi h$ in $A$.
The next input tells us what duality must do to remove that component.

\begin{theorem}[The evaluation criterion]
\label{thm:imported-evaluation}
Under the hypotheses of Theorem~\ref{thm:imported-obstruction}, suppose
that canonical biduality $I\to\mathbb D_R\mathbb D_RI$ is an
isomorphism and that canonical evaluation
\[
 I\otimes_R^{\mathbf L}\mathbb D_RI\longrightarrow R
\]
factors through $\iota:I\hookrightarrow R$ in $\mathcal D(R)$.
Then, for every $m\geq1$, $\iota_*o_m(f)=0$ implies $o_m(f)=0$.
Consequently the three invariants in \eqref{eq:three-inequalities}
are equal, including the value $\infty$.
\end{theorem}
The first criterion combines \cite[Theorem~3.3]{CV2017} with
\cite[Theorem~3.1]{ParentI}; the second is
\cite[Theorem~4.5 and Corollary~4.6]{ParentI}.
Suppose $\iota_*o_m(f)=0$. The evaluation pairing returns a dual calculation to the original
module. If its values lie in $I$, that return takes place inside $I$;
canonical biduality identifies the resulting map with the original
inclusion $I^m\to I$. This explains how a vanishing first seen in $R$
can be recovered in $I$.

The remaining inputs from \cite{ParentI} are the split-model lemma
\cite[Lemma~2.7]{ParentI}, transport of canonical duality along a
quasi-isomorphism \cite[Lemma~4.7]{ParentI}, and the general bound
$\relcat(f_\Q)\leq\msecat(f_\Q)+1$ \cite[Theorem~1.3]{ParentI}.
Our relative Frobenius argument extends the finite pairing argument
of \cite[Lemma~4.4]{ParentI}. The bounded truncation used in
\cite[proof of Corollary~5.9]{ParentI} is proved again in
Lemma~\ref{lem:bounded-truncation} below.

The higher-diagonal comparison is
\cite[Corollary~5.2]{ParentI}, recalled in the introduction and
Remark~\ref{rem:pd-fibre-diagonal}.

\section{The Thom criterion and self-intersection}
\label{sec:pd-fibres}

Fibre Poincar\'e duality gives canonical biduality for the kernel
ideal. It does not by itself make evaluation take values in that ideal.
We first state the duality input and make the Thom morphism concrete.
Its compatibility with the ideal will be measured by one complementary
component on a derived tensor product.
Throughout this section we use \eqref{eq:split-model} and assume that
$H^*(F;\Q)$ is a finite-dimensional Poincar\'e duality algebra: it has
an orientation in degree $N$ whose multiplication pairings
$H^j(F)\otimes H^{N-j}(F)\to\Q$ are nondegenerate.

\subsection{Thom morphisms and a first calculation}

The following duality statement produces the orientation morphism.
Its proof is in Appendix~\ref{app:relative-duality}.

\begin{lemma}[Duality over the base]\label{lem:relative-pd-duality}
Suppose that $H^*(F;\Q)$ is a finite-dimensional Poincar\'e duality
algebra of formal dimension $N$. Then $R$ and $I$ are perfect as
differential $A$-modules. There is an $A$-linear cocycle
$\Pi\in\operatorname{Hom}_A^{-N}(R,A)$ inducing the orientation
on the fibre, called fibre integration, whose pairing induces an
isomorphism
\begin{equation}\label{eq:relative-pd-isomorphism}
 R\xrightarrow{\simeq}
       \mathbf R\operatorname{Hom}_A(R,A)[-N]
       \qquad\text{in }\mathcal D(R).
\end{equation}
For every differential $R$-module $M$, it induces a natural
isomorphism
\begin{equation}\label{eq:relative-pd-coinduction}
 \mathbb D_RM\simeq
       \mathbf R\operatorname{Hom}_A(M,A)[-N].
\end{equation}
In particular, the natural biduality morphism
$I\to\mathbb D_R\mathbb D_R I$ is an isomorphism, and
$\mathbb D_R A\simeq A[-N]$ in $\mathcal D(R)$, with $R$ acting
on $A$ through $\varphi$.
\end{lemma}

For an $A$-module $V$, write $\varphi^*V$ for its restriction
along $\varphi:R\to A$, and use the same notation for morphisms.
A \emph{Thom morphism} is a generator $\tau$ of the
one-dimensional space
\[
 \operatorname{Hom}_{\mathcal D(R)}(A[-N],R)
       \cong H^0(A)=\Q.
\]
Its \emph{self-intersection class} is the class $e\in H^N(A)$
defined by $\iota_A^*(\varphi\tau)=m_e$ in $\mathcal D(A)$,
where $\iota_A:A\hookrightarrow R$ is the split inclusion and
$m_e:A[-N]\to A$ is multiplication by $e$. Equivalently,
$e=H^N(\varphi)(H^N(\tau)([1]))$, where
$[1]\in H^N(A[-N])$ is the shifted unit.
Replacing $\tau$ by $c\tau$, for $c\in\Q^\times$, replaces $e$ by
$ce$; the vanishing and self-intersection conditions below are
unchanged by this simultaneous rescaling.
For every $R$-module $M$, define \emph{Thom multiplication} on $M$ by
\[
 T_M=1_M\otimes\tau:
 M\otimes_R^{\mathbf L}A[-N]\longrightarrow M.
\]

\begin{example}[An even-sphere section]\label{ex:even-sphere-thom}
Let $R=A\otimes\Lambda(x_{2d},y_{4d-1})$, with $d\geq1$,
$dx=0$, $dy=x^2+ax$, and $a\in A^{2d}$ a cocycle; the map
$\varphi:R\to A$ kills $x$ and $y$.  Set $t=sx$ and $v=sy$.  Use the semifree acyclic closure $\varepsilon:\mathcal C=R\otimes\Lambda(t,v)\to A$,
with $\varepsilon|_R=\varphi$ and $\varepsilon(t)=\varepsilon(v)=0$. Define $\tau$ by
\[
 Dt=x,\qquad Dv=y-(x+a)t,\qquad
 \tau(1)=x+a,\quad\tau(t)=y,\quad
 \tau(v^j)=\tau(tv^j)=0\quad(j\geq1).
\]
Here $|\tau|=2d$ and $v$ is even, so the free basis consists of all
$v^j$ and $tv^j$. On these basis elements,
\begin{align*}
 D(v^j)&=j\bigl(y-(x+a)t\bigr)v^{j-1} &&(j\geq1),\\
 D(tv^j)&=xv^j+jytv^{j-1} &&(j\geq0),\\
 \tau D(v^j)&=j\delta_{j1}\bigl(y(x+a)-(x+a)y\bigr)=0
       =d\tau(v^j) &&(j\geq1),\\
 \tau D(tv^j)&=\delta_{j0}x(x+a)
       =\delta_{j0}dy=d\tau(tv^j) &&(j\geq0).
\end{align*}
The term multiplied by $j$ is omitted when $j=0$. In the last
line the possible term $y^2$ vanishes because $y$ is odd, and the
value on $1$ is closed. Thus $\varphi\tau=a\varepsilon$ on the entire free
basis, so its self-intersection class is $e=[a]$.
Its fibre restriction represents the orientation of the even sphere,
so this is a Thom cocycle. The derived self-intersection formula
\eqref{eq:thom-self-intersection} therefore holds with $e=[a]$,
which is nonzero precisely when $[a]\ne0$. By Theorem~\ref{thm:thom-evaluation},
$\msecat(f_\Q)=\secat(f_\Q)=\relcat(f_\Q)$ for the section represented by this model.
\end{example}

\subsection{Projection and the Thom criterion}

The Thom morphism converts restriction to the section into
multiplication. A map to $R$ becomes zero after applying $\varphi$
precisely when it vanishes after Thom multiplication. Applied to
evaluation on $I$, this is the bridge to the ideal-valued criterion.
Write $\Psi_\varphi=A\otimes_R^{\mathbf L}-$ for extension of scalars;
its right adjoint is restriction $\varphi^*$. We retain this notation
throughout the proof.

\begin{lemma}[Projection lemma]\label{lem:thom-projection}
Under the hypotheses of Lemma~\ref{lem:relative-pd-duality}, for every
$R$-module $P$ and every morphism $\beta:P\to R$ in $\mathcal D(R)$,
\begin{equation}\label{eq:thom-projection-property}
 \varphi\beta=0\quad\Longleftrightarrow\quad\beta T_P=0.
\end{equation}
\end{lemma}
\begin{proof}
The natural $A$-action on $\mathbb D_RA$ is induced by
precomposition. Compute it on a semifree $R$-module resolution
$\mathcal C_A\to A$, with its commuting $A$-action induced by $\iota_A$.
For a cocycle representing $\tau$, the map
\[
 A[-N]\longrightarrow\operatorname{Hom}_R(\mathcal C_A,R),
 \qquad a\longmapsto\tau(a\cdot-),
\]
with the standard shifted signs, is $A$-linear. For a homogeneous
$R$-linear cochain $v$, its precomposition action satisfies
$(a\cdot v)(x)=\pm v(\iota_A(a)x)=\pm\iota_A(a)v(x)$.
Thus this is also the $A$-action obtained by restriction along $\iota_A$. On cohomology it
identifies $H(A)[-N]$ with $\operatorname{Ext}_R(A,R)$:
Lemma~\ref{lem:relative-pd-duality} identifies the latter with a
free rank-one module, whose $H(R)$-action factors through
$H(\varphi)$. Its degree-$N$ part is
$\operatorname{Ext}_R^N(A,R)\cong H^0(A)=\Q$, so the nonzero
class $[\tau]$ generates this free rank-one module.
Since $\varphi\iota_A=\id_A$, the resulting $H(A)$-action is
precisely the precomposition action above.
To identify this action with coinduction, choose a homotopically injective
resolution $R\to\widehat R$. The $A$-linear zigzag
\[
 \operatorname{Hom}_R(\mathcal C_A,R)\longrightarrow
 \operatorname{Hom}_R(\mathcal C_A,\widehat R)\longleftarrow
 \operatorname{Hom}_R(A,\widehat R)
\]
consists of quasi-isomorphisms: $\mathcal C_A$ is homotopically projective and $\widehat R$
is homotopically injective. The rightmost complex has the natural
coinduction action. Thus the preceding map gives
$A[-N]\simeq\mathbf R\operatorname{Hom}_R(A,R)$ in
$\mathcal D(A)$, with adjoint $\tau$.

The right adjunction for $\varphi^*$ consequently gives a bijection
\[
 \begin{gathered}
 \operatorname{Hom}_{\mathcal D(A)}(V,A)
 \xrightarrow{\simeq}
 \operatorname{Hom}_{\mathcal D(R)}(\varphi^*(V[-N]),R),\\[2pt]
 h\longmapsto\tau \varphi^*(h[-N]).
 \end{gathered}
\]
The left adjunction sends $\varphi\beta$ to $\Psi_{\varphi}\beta:\Psi_{\varphi} P\to A$.
With $V=\Psi_{\varphi} P$, its image under the displayed bijection is
$\beta T_P$, using the graded symmetry
$\varphi^*(\Psi_{\varphi} P)[-N]\simeq P\otimes_R^{\mathbf L}A[-N]$.
This proves \eqref{eq:thom-projection-property}.
\end{proof}

\Needspace*{15\baselineskip}
\begin{theorem}[The Thom criterion]
\label{thm:thom-evaluation}
Under the hypotheses of Lemma~\ref{lem:relative-pd-duality},
choose a Thom morphism $\tau:A[-N]\to R$ with self-intersection
class $e\in H^N(A)$.
The following conditions are equivalent:
\begin{enumerate}[label=\textup{(\roman*)},leftmargin=*,itemsep=2pt]
\item Canonical derived evaluation
$I\otimes_R^{\mathbf L}\mathbb D_RI\to R$ factors through $I$.
\item $T_I=0$ in $\mathcal D(R)$.
\item The self-intersection formula holds:
\begin{equation}\label{eq:thom-self-intersection}
 \varphi\tau=\varphi^*(m_e)
       \qquad\text{in }\mathcal D(R).
\end{equation}
\end{enumerate}
Each condition implies, for every $m\geq1$,
\[
 \iota_*o_m(f)=0\quad\Longrightarrow\quad
 o_m(f)=0
 \qquad\text{in }\mathcal D(R),
\]
and consequently
\[
 \msecat(f_{\Q})=\secat(f_{\Q})=\relcat(f_{\Q}).
\]
\end{theorem}
When $e=0$, condition \textup{(iii)} is $\varphi\tau=0$ in
$\mathcal D(R)$, an identity before restriction along $\iota_A$.
\begin{proof}
\emph{Equivalence of \textup{(i)} and \textup{(ii)}.}
Put $D_I=\mathbb D_RI$ and write
$\operatorname{ev}:I\otimes_R^{\mathbf L}D_I\to R$.
If $T_I=0$, naturality and graded symmetry give
$T_{I\otimes_R^{\mathbf L}D_I}=T_I\otimes1_{D_I}=0$.
Lemma~\ref{lem:thom-projection} gives $\varphi\operatorname{ev}=0$;
the exact triangle $I\to R\to A$ then factors evaluation through
$I$. Conversely, such a factorization gives
$\varphi\operatorname{ev}=0$ and hence
$\operatorname{ev}(T_I\otimes1_{D_I})=0$ by the same property.
Its tensor--Hom adjoint is $\eta_I T_I$.
Lemma~\ref{lem:relative-pd-duality} makes $\eta_I$ an isomorphism,
so $T_I=0$.

\emph{Use of the two tensor units.}
The self-intersection class determines the restriction to the base,
but the identity to prove is $R$-linear. Under adjunction, its extra
content is the vanishing of a complementary component of Thom
multiplication. The following two units keep these linearities distinct.
Set
\[
 \Gamma=A\otimes_R^{\mathbf L}A,\qquad
 \mu:\Gamma\longrightarrow A.
\]
There are two unit maps, with different categories of linearity:
\begin{equation}\label{eq:two-tensor-units}
 \begin{aligned}
 \varphi\otimes1_A &: \varphi^*A\longrightarrow \varphi^*\Gamma
          &&\text{in }\mathcal D(R),\\
 1_A\otimes\varphi &: A\longrightarrow\Gamma
          &&\text{in }\mathcal D(A).
 \end{aligned}
\end{equation}
The first is the unit of extension--restriction adjunction at $\varphi^*A$.
The second is $\Psi_{\varphi}(\varphi)$. The two splitting identities are
\begin{equation}\label{eq:two-unit-splittings}
 \varphi^*(\mu)(\varphi\otimes1_A)=1_{\varphi^*A},
 \qquad \mu(1_A\otimes\varphi)=1_A.
\end{equation} Applying $\Psi_{\varphi}$ to $I\to R\to A$ therefore gives
\[
 \Gamma\simeq A\oplus\mathcal Z,
 \qquad \mathcal Z=\Psi_{\varphi}(I)[1].
\]
The $A$-linear Thom multiplication
$\widetilde T_A=\Psi_{\varphi}(\tau):\Gamma[-N]\to A$ is adjoint to
$\varphi\tau$ and satisfies $\varphi^*(\widetilde T_A)=T_{\varphi^*A}$.
Write the triangle of the augmentation as
\[
 I\xrightarrow{\iota}R\xrightarrow{\varphi}\varphi^*A
       \xrightarrow{\delta}I[1].
\]
Tensor this triangle on the right with $A[-N]$.
Naturality of Thom multiplication along the $R$-linear boundary
$\delta$ gives, up to sign,
\begin{equation}\label{eq:thom-boundary-naturality}
 \delta \varphi^*(\widetilde T_A)=T_I[1]b_1,
 \qquad b_1:\varphi^*(\Gamma[-N])\to I\otimes_R^{\mathbf L}A[1-N].
\end{equation}
The middle arrow of this tensor triangle is
$(\varphi\otimes1_A)[-N]$.
Since $\varphi^*(\mu[-N])$ retracts it, $b_1$ is a split epimorphism.

\emph{From \textup{(iii)} to \textup{(ii)}.}
Adjunction sends \eqref{eq:thom-self-intersection} to
$\widetilde T_A=m_e\mu[-N]$ in $\mathcal D(A)$.
Since $\delta \varphi^*(m_e)=\delta\varphi\tau=0$, it follows that
$\delta \varphi^*(\widetilde T_A)=0$.
Equation~\eqref{eq:thom-boundary-naturality} and a right inverse
to $b_1$ give $T_I=0$.

\emph{From \textup{(ii)} to \textup{(iii)}.}
The splitting $\Gamma\simeq A\oplus\mathcal Z$ reduces the
assertion to the vanishing of the complementary component of
$\widetilde T_A$. Restriction along $\iota_A$ will then identify
the remaining component with $m_e$.
Apply $\Psi_{\varphi}$ to the augmentation triangle. The resulting
triangle in $\mathcal D(A)$ is
\[
 \Psi_{\varphi}(I)\longrightarrow
 A\xrightarrow{1_A\otimes\varphi}\Gamma
       \xrightarrow{b_2}\mathcal Z.
\]
Because $\mu(1_A\otimes\varphi)=1_A$, it splits. Choose a
section $j_0$ of $b_2$ and put
$j=j_0-(1_A\otimes\varphi)\mu j_0$. Then
\[
 b_2j=1_{\mathcal Z},\qquad \mu j=0,\qquad
 1_\Gamma=(1_A\otimes\varphi)\mu+jb_2.
\]
All these identities hold in $\mathcal D(A)$.
If $T_I=0$, equation~\eqref{eq:thom-boundary-naturality} gives
$\delta \varphi^*(\widetilde T_A)=0$. The complementary component
\[
 \kappa_{\mathrm{Th}}=\widetilde T_Aj[-N]:\mathcal Z[-N]\longrightarrow A
       \qquad\text{in }\mathcal D(A)
\]
therefore satisfies $\delta \varphi^*\kappa_{\mathrm{Th}}=0$.
Exactness of $\operatorname{Hom}_{\mathcal D(R)}(\varphi^*(\mathcal Z[-N]),-)$
gives a lift
\[
 \gamma:\varphi^*(\mathcal Z[-N])\longrightarrow R,
 \qquad \varphi\gamma=\varphi^*\kappa_{\mathrm{Th}}.
\]
Associativity of the derived tensor product and the shift
identifications give
\[
 T_{\varphi^*(\mathcal Z)}\simeq \varphi^*(\Psi_{\varphi}(T_I))[1]=0.
\]
Before shifting, both morphisms are represented by
$1_A\otimes1_I\otimes\tau$ from
$A\otimes_R^{\mathbf L}I\otimes_R^{\mathbf L}A[-N]$
to $A\otimes_R^{\mathbf L}I$.
Hence $T_{\varphi^*(\mathcal Z[-N])}=0$, and the projection property
applied to $\gamma$ yields $\varphi^*\kappa_{\mathrm{Th}}=0$.
Restriction along $\iota_A$ is a left inverse to $\varphi^*$, since
$\varphi\iota_A=\id_A$. Thus $\kappa_{\mathrm{Th}}=0$ in $\mathcal D(A)$ itself.
The splitting now gives
\[
 \widetilde T_A=(\widetilde T_A(1_A\otimes\varphi)[-N])\mu[-N]
       \qquad\text{in }\mathcal D(A).
\]
Let $e'\in H^N(A)$ be the class represented by
$\widetilde T_A(1_A\otimes\varphi)[-N]$. The two sides are the adjoints of
$\varphi\tau$ and $\varphi^*(m_{e'})$, respectively, so
$\varphi\tau=\varphi^*(m_{e'})$. Restriction along $\iota_A$ identifies
$m_{e'}=\iota_A^*(\varphi\tau)=m_e$.
This proves \eqref{eq:thom-self-intersection} with the prescribed
class $e$.

Finally, Lemma~\ref{lem:relative-pd-duality} supplies canonical
biduality, and condition \textup{(i)} supplies ideal-valued evaluation.
Theorem~\ref{thm:imported-evaluation} gives the assertion about ideal
powers and the equality of the three invariants.
\end{proof}

The validity of these equivalent conditions for an arbitrary
Poincar\'e duality fibre remains open.
We next prove them under a degree bound. Later sections obtain them
from finite total rational homotopy, smooth tubular neighbourhoods,
and relative Frobenius pairings.

\subsection{A dimension bound that kills the obstruction}

The proof locates $T_I$ in a derived Hom group. A lower bound on
the degrees of its source and an upper bound on the cohomology of
its target force that group to vanish.

The following truncation turns a cohomological bound into a bound on
cochains (it is also used in \cite[proof of Corollary~5.9]{ParentI}).

\begin{lemma}[Bounded truncation]\label{lem:bounded-truncation}
Let $A$ be a nonnegatively graded CDGA with $A^0=\Q$, and suppose
that $H^{>n}(A)=0$ for an integer $n\geq0$. Choose a vector-space
complement $L'$ to $\ker(d:A^n\to A^{n+1})$ in $A^n$, taking
$L'=0$ when $n=0$. Then
\[
 \overline A=A/(L'\oplus A^{>n})
\]
is a CDGA concentrated in degrees at most $n$, and the quotient
$A\to\overline A$ is a quasi-isomorphism.
If $R\to A$ is a CDGA morphism and $P$ is a semifree $R$-module
whose cells have degrees at least $c$, then
\[
 \operatorname{Hom}_R(P,\overline A)^k=0
 \qquad\text{whenever }c+k>n.
\]
Moreover, $\operatorname{Hom}_R(P,A)\to
\operatorname{Hom}_R(P,\overline A)$ is a quasi-isomorphism.
\end{lemma}
\begin{proof}
Connectedness and nonnegative grading make
$L'\oplus A^{>n}$ an ideal, and its differential stays in that ideal.
It has no cycles in degree $n$ because $L'$ complements the cycles.
In degree $n+1$, the identity $dL'=dA^n$ identifies its cohomology
with $H^{n+1}(A)=0$; in higher degrees its cohomology equals that
of $A$ and is zero. The ideal is therefore acyclic, proving the
first assertion. A homogeneous $R$-linear map of degree $k$ sends
a cell of degree $p\geq c$ into $\overline A^{p+k}=0$ if
$c+k>n$, which proves the vanishing. The last assertion follows
because Hom from a semifree module preserves quasi-isomorphisms.
\end{proof}

\begin{theorem}[Equality under a base-dimension bound]\label{thm:pd-fibre-dimension}
Let $f:X\to Y$ be a rationally retractive map between simply
connected spaces of finite rational type, and let $g:Y_{\Q}\to
X_{\Q}$ be a homotopy retraction. Suppose that the rational
cohomology of the homotopy fibre of $g$ is a finite-dimensional
Poincar\'e duality algebra of formal dimension $N$.
In the associated split minimal model, suppose $W^{<c}=0$ for
an integer $c\geq2$. If there is an integer $n\geq0$ such that
\[
 H^{>n}(X;\Q)=0,
 \qquad n<N+c-1,
\]
then
\[
 \msecat(f_{\Q})=\secat(f_{\Q})=\relcat(f_{\Q}).
\]
In particular, it suffices that the cohomological dimension of
$X$ be at most $N$.
\end{theorem}
\begin{proof}
\emph{Step 1: locate the possible obstruction.}
Choose a Thom morphism $\tau:A[-N]\to R$.
The map $I\otimes_R^{\mathbf L}A\to A$ induced by
$I\hookrightarrow R$ is zero. Indeed, its exact triangle continues
with the $R$-linear unit
$\varphi\otimes1_A:\varphi^*A\to \varphi^*\Gamma$, whose left inverse is
$\varphi^*(\mu)$ by \eqref{eq:two-unit-splittings}; hence the first map is zero.
Naturality therefore gives
\[
 \iota T_I=\tau(\iota\otimes1_{A[-N]})=0,
 \qquad \iota:I\hookrightarrow R.
\]
Apply $\operatorname{Hom}_{\mathcal D(R)}
(I\otimes_R^{\mathbf L}A[-N],-)$ to the exact triangle
$I\to R\to A$. To prove $T_I=0$, it suffices to show that
\begin{equation}\label{eq:pd-fibre-factor-group}
 \operatorname{Hom}_{\mathcal D(R)}
   (I\otimes_R^{\mathbf L}A[-N],A[-1])
 \cong
 \operatorname{Hom}_{\mathcal D(R)}
   (A\otimes_R^{\mathbf L}I,A[N-1])
 =0.
\end{equation}

\emph{Step 2: use the degree bounds.}
Since $I$ is concentrated in degrees at least $c$, it has a
semifree resolution $P_1\to I$ whose cells have degrees at least $c$.
The construction is the standard cell-adjoining resolution
\cite[Proposition~6.6 and \S6, Exercise~8]{FHT2001}; the degree bound follows as below.
In the standard inductive construction, start with representatives
of cohomology in degrees at least $c$ and kill independent kernel
classes. The degree-$c$ killing cells have independent differentials
in degree $c+1$, so they create no new degree-$c$ cycles; hence no
cell below degree $c$ is needed.
The cone
$Q_A=R\oplus sP_1$, with
$D(sp)=\varepsilon_1(p)-sDp$, is a semifree resolution of $A$.
Its cells have degree zero or at least $c-1$.
Thus $Q_A\otimes_RP_1$ models
$A\otimes_R^{\mathbf L}I$ with every cell in degree at least $c$.
Apply Lemma~\ref{lem:bounded-truncation} to $A$, the source
$Q_A\otimes_RP_1$, and the degree $N-1$. It gives a
quasi-isomorphic target $\overline A$ concentrated in degrees at
most $n$, for which
\[
 \operatorname{Hom}_R(Q_A\otimes_RP_1,\overline A)^{N-1}=0
 \qquad\text{since }c+N-1>n.
\]
This proves \eqref{eq:pd-fibre-factor-group}, so $T_I=0$.
Theorem~\ref{thm:thom-evaluation} gives the equality of invariants.
Finally, $c\geq2$ makes $n\leq N$ a sufficient dimension bound.
\end{proof}

\begin{remark}[A second proof for Poincar\'e duality diagonals]
\label{rem:pd-fibre-diagonal}
For the higher diagonal $X\to X^r$, $r\geq2$, the first projection is
a retraction with fibre $X^{r-1}$. If $H(X;\Q)$ has Poincar\'e
formal dimension $n$, then the fibre has dimension $N=(r-1)n$.
Taking $c=2$ gives $n<N+1$ for every $n\geq0$ and $r\geq2$, so
Theorem~\ref{thm:pd-fibre-dimension} proves equality of all three
invariants of the diagonal. Since
$\TC_r(X_\Q)=\secat((\Delta_r)_\Q)$ and
$\MTC_r(X)=\msecat((\Delta_r)_\Q)$, this recovers
$\MTC_r(X)=\TC_r(X_\Q)=\relcat((\Delta_r)_\Q)$,
the higher-diagonal equality of \cite[Corollary~5.2]{ParentI}.
For a general fibre outside the dimension range, the group in
\eqref{eq:pd-fibre-factor-group} can be nonzero.
\end{remark}

\section{Equality for fibres with finite total rational homotopy}
\label{sec:finite-homotopy}

Finite total rational homotopy gives another way to verify the
Thom criterion. No bound on the cohomological dimension of the base
is required.

\begin{theorem}[Fibres with finite total rational homotopy]
\label{thm:finite-homotopy-fibre}
Let $f:X\to Y$ be a rationally retractive map between simply
connected spaces of finite rational type, and choose a homotopy
retraction $g:Y_{\Q}\to X_{\Q}$. If its homotopy fibre $F$ satisfies
\[
 \dim_{\Q}\bigl(\pi_*(F)\otimes\Q\bigr)<\infty,
\]
then
\[
 \msecat(f_{\Q})=\secat(f_{\Q})=\relcat(f_{\Q}).
\]
The equality includes the value $\infty$ and requires no bound on
the cohomological dimension of $X$.
\end{theorem}

The proof separates three assertions. First, a finite module threshold
forces the fibre to have finite cohomology. Second, a finite relative
basis gives a generator of the relative Ext module. Third, a length
estimate makes the augmentation of this generator scalar multiplication.
The first assertion supplies the elliptic inequality
$n_{\mathrm{odd}}\geq n_{\mathrm{even}}$ needed by the third.

Reduction modulo $A^+$ leaves no scalar-linear term in the relative
differential, so $(\Lambda W,\overline d)$ is minimal. For the fibre
of $g$, the homotopy interpretation of a minimal model gives
\begin{equation}\label{eq:geometric-W}
 W^j\cong\operatorname{Hom}_{\Q}(\pi_j(F)\otimes\Q,\Q),
 \qquad j\geq2;
\end{equation}
see \cite[Proposition~15.5 and Theorem~15.11]{FHT2001}.
The hypothesis is therefore $\dim W<\infty$, not merely finite
dimension in each degree. The split homotopy sequence identifies
$\pi_*(Y_\Q)$ with $\pi_*(X_\Q)\oplus\pi_*(F)$; only the added
summand is required to have finite total dimension.

\subsection{First stage: a finite threshold forces ellipticity}

Finitely many homotopy generators can still yield infinitely many
cohomology classes. We first use the module nullhomotopy to rule out
that possibility at a finite threshold.

\begin{lemma}\label{lem:finite-threshold-elliptic}
If $\dim W<\infty$ and $\msecat(f_\Q)<\infty$, then
$H^*(\Lambda W,\overline d)$ is finite dimensional. Thus $F$ is
rationally elliptic.
\end{lemma}
\begin{proof}
We must retain the module nullhomotopy after passing to the fibre.
Choose $m\geq1$ such that
$[I^m\hookrightarrow R]=0$ in $\mathcal D(R)$, using
Theorem~\ref{thm:imported-obstruction}. The powers $I^m$ are semifree
over $A$. To see this, choose a finite Sullivan basis
$w_1,\ldots,w_k$ and assign positive integer weights recursively so
that the weight of $w_i$ exceeds the weight of each relative monomial
in $d_Rw_i$. Apart from the differential of the $A$-coefficient,
the differential strictly lowers this weight and preserves relative
length at least $m$. Increasing weight orders a semifree relative
monomial basis of $I^m$.

Choose a semifree $R$-resolution $P_m\to I^m$. It is also semifree
over $A$, and the displayed quasi-isomorphism is an $A$-module
homotopy equivalence. Tensoring it and the module nullhomotopy over
$A\to\Q$ therefore gives
\[
 [\Lambda^{\geq m}W\hookrightarrow(\Lambda W,\overline d)]=0
       \quad\text{in }\mathcal D(\Lambda W,\overline d).
\]
If $W=0$ there is nothing to prove. Otherwise let $D$ be the largest
generator degree. Every cochain of degree greater than $(m-1)D$ has
length at least $m$. A cocycle in such a degree is consequently in
$\Lambda^{\geq m}W$, and the displayed nullity makes its class in
$H(\Lambda W,\overline d)$ zero. Hence $H^{>(m-1)D}(\Lambda W,\overline d)=0$. Each remaining degree is finite
dimensional, which proves the assertion.
\end{proof}

\subsection{Second stage: a relative Thom cocycle}

A Thom cocycle represents the orientation generator in relative Ext.
We construct one by adjoining the relative Sullivan generators
successively, refining the diagonal construction of
\cite[Definition~5.1 and Lemma~5.2]{Wakatsuki2016} to respect the
augmentation and relative length. Relative minimality excludes
scalar-linear terms in their differentials; the comparison details
are in Appendix~\ref{app:filtered-thom}.

\begin{lemma}[A relative Thom cocycle]\label{lem:filtered-thom}
Suppose that $W$ is finite dimensional.  Write
$n_{\mathrm{odd}}=\dim W^{\mathrm{odd}}$,
$n_{\mathrm{even}}=\dim W^{\mathrm{even}}$, and put
\[
 N=\sum_{w\in\mathcal W_{\mathrm{odd}}}|w|
       -\sum_{w\in\mathcal W_{\mathrm{even}}}(|w|-1),
\]
where $\mathcal W$ is a homogeneous basis of $W$.
There is a relative acyclic closure, that is, a relative Sullivan
extension resolving the augmentation as in \cite[Sections~14--15]{FHT2001},
$\varepsilon:\mathcal C_A=R\otimes\Lambda(sW)\longrightarrow A$ and a cocycle
$\tau\in\operatorname{Hom}_R^N(\mathcal C_A,R)$ such that
\[
 H(A)\,[\tau]\cong \operatorname{Ext}_R(A,R).
\]
\end{lemma}
Here $N$ is the degree of the Ext generator. Without ellipticity it
may be negative: for one even generator of degree $2d$, it is $1-2d$.
When the fibre is elliptic it equals its formal dimension.

\begin{proof}
\emph{Step 1: the acyclic closure.}
Order the basis $w_1,\ldots,w_k$ by nondecreasing degree.
This is a Sullivan order: minimality and $A^1=0$ imply that every
relative generator occurring in $d_Rw_i$ has degree less than $|w_i|$.
Set
$R_{<i}=A\otimes\Lambda(w_1,\ldots,w_{i-1})$.
Lemma~\ref{lem:relative-thom-cylinder} supplies an acyclic closure with
\[
 D(sw_i)=w_i-z_i,\qquad Dz_i=dw_i,\qquad
 z_i\in R_{<i}^+\otimes\Lambda(sw_1,\ldots,sw_{i-1}),
\]
in which every monomial of $z_i$ contains a letter from $W\oplus sW$.
The construction is
compatible with every initial segment of the Sullivan basis: for
$R_i=A\otimes\Lambda(w_1,\ldots,w_i)$, the same terms $z_j$,
$j\leq i$, make
$\mathcal C_i=R_i\otimes\Lambda(sw_1,\ldots,sw_i)$ an acyclic closure of
$A$ over $R_i$. Thus the induction uses compatible resolutions
at each stage.

\emph{Step 2: adjoining one relative generator.}
Suppose that the closure $\mathcal C_{<i}$ is already constructed over
${R_{<i}}=A\otimes\Lambda(w_1,\ldots,w_{i-1})$, and write
$x=w_i$, $t=sx$, $R_i={R_{<i}}\otimes\Lambda x$, and $Dt=x-z$.
Extend a homogeneous ${R_{<i}}$-linear cochain $\xi:\mathcal C_{<i}\to {R_{<i}}$ of degree $p$
to $R_i\otimes_{R_{<i}}\mathcal C_{<i}$ by graded $R_i$-linearity.
For homogeneous monomials $\omega$ in the old suspended generators, define
$\Phi(\xi)$ by graded $R_i$-linearity and the following formulas:
\begin{align}\label{eq:thom-phi-odd}
 |x|\text{ odd}:&\qquad
 \Phi(\xi)(\omega)=x\xi (\omega)-(-1)^p \xi(z\omega ),&
 \Phi(\xi)(\omega t^j)&=0\quad(j\geq1),\\
 \label{eq:thom-phi-even}
 |x|\text{ even}:&\qquad
 \Phi(\xi)(\omega)=0,&
 \Phi(\xi)(\omega t)&=(-1)^{p+|\omega|}\xi(\omega).
\end{align}
The operator $\Phi$ has degree $|x|$ in the odd case and $1-|x|$
in the even case. Both degrees are odd, and
$\partial\Phi=-\Phi\partial$.
Lemma~\ref{lem:complete-thom-comparison} proves that this is a
quasi-isomorphism after the indicated shift, using the filtration by
the degree of the ${R_{<i}}$-coefficients.  Starting with the identity cochain
of $A$ and applying these comparisons therefore gives a cocycle
$\tau$ generating $\operatorname{Ext}_R(A,R)$ over $H(A)$, of degree $N$.
\end{proof}

\subsection{Third stage: the length estimate and self-intersection}

The cocycle must do more than generate Ext: it must vanish after
augmentation on every suspended monomial except the unit. The two
parities in the construction affect the length estimate in opposite
ways, and ellipticity gives exactly the inequality that is needed.

\begin{lemma}[Relative-length estimate]\label{lem:thom-length}
The cocycle in Lemma~\ref{lem:filtered-thom} satisfies
\[
 \tau(\omega)\in I^{\ell(\omega)+n_{\mathrm{odd}}-n_{\mathrm{even}}}
 \qquad(\omega\text{ a monomial in }\Lambda(sW)),
\]
where $\ell(\omega)$ is word length and $I^j=R$ for $j\leq0$.
\end{lemma}
\begin{proof}
Suppose at one stage that
$\xi(\omega)\in I_{<i}^{\ell(\omega)+n_{\mathrm{odd},<i}-n_{\mathrm{even},<i}}$, and put $c=n_{\mathrm{odd},<i}-n_{\mathrm{even},<i}$.
For odd $x$, the term $x\xi (\omega)$ has the required extra relative factor.
A term of $z$ with $a$ unsuspended relative letters and $b$ suspended
letters satisfies $a+b\geq1$.  Its contribution to $\xi(z\omega )$ lies in
\[
 I^{a+\max\{0,\ell(\omega)+b+c\}}\subseteq I^{\ell(\omega)+c+1}
\]
when the exponent on the right is positive. Indeed,
$\ell(\omega)+c+1>0$ implies $\ell(\omega)+c\geq0$, so the exponent on
the left is $\ell(\omega)+c+a+b\geq\ell(\omega)+c+1$.
For a nonpositive exponent the assertion is automatic.
For even $x$, the only nonzero new value is, up to sign, $\xi(\omega)$ on
$\omega t$; its input length increases by one while
$n_{\mathrm{odd}}-n_{\mathrm{even}}$ decreases by one.
This proves the estimate by induction, including the convention
$I^j=R$ for $j\leq0$.
\end{proof}

\begin{corollary}\label{cor:finite-thom-self-intersection}
For the cocycle of Lemma~\ref{lem:thom-length}, suppose that
$n_{\mathrm{odd}}\geq n_{\mathrm{even}}$.  Then the cocycle can be chosen so that
\[
 \varphi\tau=(e_A\cdot-)\varepsilon,
 \qquad e_A=\varphi(\tau(1))\in A^N,
 \qquad de_A=0.
\]
In particular, the self-intersection formula holds already at the
level of cochains. Under the Poincar\'e duality hypothesis of
Theorem~\ref{thm:thom-evaluation}, the self-intersection class
of that theorem is $e=[e_A]$.
If $n_{\mathrm{odd}}>n_{\mathrm{even}}$, then $e_A=0$.
\end{corollary}

\begin{proof}
The ideal-power estimate puts $\tau(\omega)$ in $I$ whenever $\omega\ne1$.
On $1$ the formula is the definition of $e_A$; graded linearity
extends it to all of $\mathcal C_A$. Since $\tau$ is a cocycle and $D1=0$,
$e_A$ is closed. Under the fibre Poincar\'e duality hypothesis,
$[\tau]$ generates the one-dimensional group
$\operatorname{Ext}_R^N(A,R)\cong H^0(A)=\Q$, so it represents a
Thom morphism as defined in Section~\ref{sec:pd-fibres}.
If $n_{\mathrm{odd}}>n_{\mathrm{even}}$, the same estimate also gives
$\tau(1)\in I$.
\end{proof}

\subsection{Proof of the equality theorem}\label{subsec:finite-equality-proof}

\begin{proof}[Proof of Theorem~\ref{thm:finite-homotopy-fibre}]
If $\msecat(f_\Q)=\infty$, the lower bounds
\eqref{eq:three-inequalities} give the equality. Otherwise
Lemma~\ref{lem:finite-threshold-elliptic} makes the fibre elliptic.
An elliptic minimal algebra has Poincar\'e duality cohomology,
$n_{\mathrm{odd}}\geq n_{\mathrm{even}}$, and formal dimension
\[
 N=\sum_{|w_i|\text{ odd}}|w_i|
       -\sum_{|w_i|\text{ even}}(|w_i|-1);
\]
see \cite[Theorem~32.6, and Proposition~38.3]{FHT2001}.
Lemma~\ref{lem:filtered-thom} constructs a generator of relative Ext
in this same formal dimension. Since
$\operatorname{Ext}_R^N(A,R)\cong H^0(A)=\Q$, its nonzero class is
a Thom morphism in the sense of Section~\ref{sec:pd-fibres}.
Lemma~\ref{lem:thom-length} and
Corollary~\ref{cor:finite-thom-self-intersection} give its
self-intersection formula. Theorem~\ref{thm:thom-evaluation} now
identifies all three invariants.
\end{proof}

\subsection{An infinite value and homogeneous fibres}

\begin{example}\label{ex:homotopy-finite-infinite-value}
The projection $X\times K(\Q,2d)\to X$, for $d\geq1$,
has fibre with one rational homotopy generator and cohomology
$\Q[x]$, $|x|=2d$. Every power of the kernel class $x$ is nonzero,
so all three invariants of its canonical section are infinite. This is the
nonelliptic case of Theorem~\ref{thm:finite-homotopy-fibre}.
\end{example}

Compact homogeneous spaces have finite total rational homotopy, so Theorem~\ref{thm:finite-homotopy-fibre} immediately treats sections
of bundles with such fibres. For groups themselves one can also
read the exact value from their odd generators.

\begin{corollary}[A finite odd relative basis]
\label{cor:finite-odd-fibre}
Let $f:X\to Y$ be rationally retractive, with simply connected
spaces of finite rational type. If its split minimal model
$R=A\otimes\Lambda W\to A$ has $W=\langle w_1,\ldots,w_r\rangle$
concentrated in odd degrees, then
\[
 \msecat(f_{\Q})=\secat(f_{\Q})=\relcat(f_{\Q})=r.
\]
\end{corollary}
\begin{proof}
The case $r=0$ is a rational equivalence. For $r\geq1$, set
$I=(W)$ and $\omega=w_1\cdots w_r$.
The section and Sullivan order imply that every monomial in
$d_Rw_i$ contains a relative generator different from $w_i$.
Consequently $d_R\omega=0$. Differentiating a proper exterior
monomial cannot produce $\omega$, since its differentiated
generator disappears. The coefficient of $\omega$ in a boundary
is therefore a boundary in $A$. As $1$ is not a boundary,
$[\omega]\ne0$ in $H(R)$.
The module criterion of Theorem~\ref{thm:imported-obstruction}
gives a lower bound $r$, while $I^{r+1}=0$ and its relative
criterion give the upper bound $r$.
\end{proof}

\begin{corollary}[Homogeneous and Grassmann bundles]
\label{cor:homogeneous-fibres}
Let $g:Y\to X$ be a fibration with a section $f:X\to Y$, where $X$ and
its fibre $F$ are simply connected and of finite rational type. If $F$ is compact homogeneous, then
\[
 \msecat(f_{\Q})=\secat(f_{\Q})=\relcat(f_{\Q}).
\]
In particular, this applies to the canonical section of
\[
 \operatorname{Gr}_a(\underline{\mathbb C}^{\,a}\oplus\eta)
 \longrightarrow X
\]
for a rank-$b$ complex vector bundle $\eta$ and $a,b\geq1$.
It also applies to the canonical section of $P\times_H(G/H)\to X$
for a principal $H$-bundle, a compact simply connected Lie group
$G$, and a closed connected subgroup $H$. For the identity section
of an adjoint bundle $Q\times_GG\to X$, the common value is the
number of odd rational homotopy generators of $G$.
\end{corollary}
\begin{proof}
A simply connected compact homogeneous space is rationally elliptic
\cite[Section~32\textup{(e)}, Example~2]{FHT2001}, so
Theorem~\ref{thm:finite-homotopy-fibre} applies.
The Grassmann bundle has fibre
$U(a+b)/(U(a)\times U(b))$; its trivial rank-$a$ summand gives the
section. In the associated bundle, left translation by $H$ fixes
the coset $H$, giving the section $x\mapsto[u,H]$.
Conjugation fixes the identity of $G$, so the adjoint bundle has
its identity section. A compact simply connected Lie group has a
finite odd minimal basis. By \eqref{eq:geometric-W}, the relative
generating space of the adjoint bundle is therefore odd and has
dimension $\operatorname{rank}G$; Corollary~\ref{cor:finite-odd-fibre}
gives its value. No restriction on the characteristic classes or
cohomological dimension of the base is needed.
\end{proof}

\section{Smooth sections and relative Frobenius models}
\label{sec:smooth-frobenius}

The Thom criterion leaves one concrete task: remove the complementary
component $\kappa_{\mathrm{Th}}$. A tubular neighbourhood does so geometrically,
by the self-intersection formula. A perfect multiplication pairing
does so algebraically, by forcing evaluation into the kernel. We
establish these two constructions before choosing examples.

\subsection{Sections of smooth bundles}

A smooth section has a normal bundle. Its Thom class supplies the
self-intersection formula, so the algebraic criterion applies to
closed manifold fibres even when their total rational homotopy is
infinite dimensional. The construction below also allows a
noncompact base.

\begin{theorem}[Smooth sections]\label{thm:smooth-section}
Let $g:Y\to X$ be a smooth fibre bundle over a simply connected smooth
manifold without boundary of finite rational type. Suppose its fibre is a closed simply
connected manifold and $g$ admits a smooth section $f:X\to Y$. Then
\[
 \msecat(f_{\Q})=\secat(f_{\Q})=\relcat(f_{\Q}).
\]
\end{theorem}
\begin{proof}
The homotopy sequence of the bundle shows that $Y$ is simply
connected. Since $F$ is closed, its rational cohomology is finite
dimensional; the Serre spectral sequence then makes $H^j(Y;\Q)$
finite dimensional in every degree. Thus $Y$ satisfies the standing
finite-rational-type hypothesis. The proof has three steps:
construct the Thom morphism near the section, compute its restriction,
and transport that identity to the split Sullivan model.

\emph{Step 1: the Thom morphism.}
Write $n=\dim F$ and identify $X$ with $f(X)$. The normal bundle
of this closed embedding is $\nu=f^*T_gY$, where $T_gY$ is the
vertical tangent bundle. It is orientable because
$w_1(\nu)\in H^1(X;\mathbb Z/2)=0$.
For compact oriented embeddings, the derived Thom construction
below appears in F\'elix--Thomas \cite{FelixThomas2009}.
Here a tubular neighbourhood and relative duality for the closed
fibre give the construction without a compactness assumption on
$X$.

Choose an open tubular neighbourhood $U_\nu$ of $X$ and put
\[
 \mathcal P_Y=A_{PL}(Y),\qquad
 \mathcal P_\nu=A_{PL}(U_\nu),\qquad
 \mathcal P_X=A_{PL}(X).
\]
Restriction factors $\alpha=A_{PL}(f)$ as
$\mathcal P_Y\xrightarrow{\rho}\mathcal P_\nu
\xrightarrow{a}\mathcal P_X$, with $a$ a quasi-isomorphism.
For a pair $(Z,Z')$, write $A_{PL}(Z,Z')$ for the kernel of
restriction to $Z'$. Set
\[
 K_Y=A_{PL}(Y,Y\setminus X),\qquad
 K_\nu=A_{PL}(U_\nu,U_\nu\setminus X).
\]
The relative-to-absolute maps $b_Y:K_Y\to\mathcal P_Y$ and
$b_\nu:K_\nu\to\mathcal P_\nu$ fit into a commutative square
with $\rho$ and the excision quasi-isomorphism
$q:K_Y\to K_\nu$.

Choose a cocycle $\theta\in K_\nu^n$ representing a nonzero
Thom class of $\nu$. Graded multiplication by $\theta$ gives a
quasi-isomorphism of $\mathcal P_\nu$-modules
$t_\theta:\mathcal P_\nu[-n]\to K_\nu$, by the Thom isomorphism
for an oriented vector bundle
\cite[Theorem~10.4]{MilnorStasheff1974}. The zigzag
\[
 \tau:\mathcal P_X[-n]
 \xleftarrow[\simeq]{a[-n]}\mathcal P_\nu[-n]
 \xrightarrow[\simeq]{t_\theta}K_\nu
 \xleftarrow[\simeq]{q}K_Y
 \xrightarrow{b_Y}\mathcal P_Y
\]
defines a morphism in $\mathcal D(\mathcal P_Y)$. The class
$H^n(\tau)([1])$, where $[1]\in H^n(\mathcal P_X[-n])$, restricts
to a nonzero orientation class on a fibre, so $\tau\ne0$.
After the model comparison in Step~3, relative duality will identify
this nonzero morphism as a generator of the one-dimensional Thom space.

\emph{Step 2: the self-intersection formula.}
Put $c_\nu=a(b_\nu\theta)\in\mathcal P_X^n$. This cocycle
represents the Euler class for the chosen Thom class.
Multiplicativity of $a$ gives the cochain identity
\[
 a\,b_\nu\,t_\theta
   =(c_\nu\cdot-)\,a[-n].
\]
Consequently the commutative excision square gives
\[
 \alpha\tau
 =a\,b_\nu\,t_\theta\,(a[-n])^{-1}
 =\alpha^*([c_\nu]\cdot-)
 \quad\text{in }\mathcal D(\mathcal P_Y).
\]

\emph{Step 3: passage to the algebraic criterion.}
Choose a split Sullivan model
$\varphi:R\to A$ for the section. The restriction
$\alpha:\mathcal P_Y\to\mathcal P_X$ is surjective by the extension
property of polynomial forms. Choose model maps $m_Y$ and $m_X$.
The relative Sullivan homotopy lifting property
\cite[Section~14\textup{(a)}]{FHT2001} replaces $m_Y$ by a homotopic map
for which $\alpha m_Y=m_X\varphi$. Indeed, the path-lifting map has
kernel $\ker\alpha\otimes\ker(\operatorname{ev}_0)$, which is contracted
by polynomial integration; cofibrancy supplies the lift. We obtain
\[
\begin{tikzcd}[column sep=3.5em,row sep=1.6em]
 R \arrow[r,"\varphi"] \arrow[d,"\simeq"'] & A \arrow[d,"\simeq"] \\
 \mathcal P_Y \arrow[r,"\alpha"] & \mathcal P_X\rlap{.}
\end{tikzcd}
\]
Restrict the identity from Step~2 along $m_Y$. The strict equality
$\alpha m_Y=m_X\varphi$ makes $m_X$ and $m_Y$ quasi-isomorphisms
of the relevant $R$-modules. Inverting these maps in $\mathcal D(R)$
therefore transports both the Thom morphism and its identity to the
split model. Since $F$ is closed, its cohomology has Poincar\'e duality. The
space of morphisms $A[-n]\to R$ is consequently one dimensional
by Lemma~\ref{lem:relative-pd-duality}. The transported morphism
is nonzero by its fibre orientation in Step~1, so it is a Thom
morphism. The transported formula is
condition~\textup{(iii)} of Theorem~\ref{thm:thom-evaluation},
with $e$ the Euler class of $\nu$; the theorem gives the
equality, including the infinite case.
\end{proof}

\subsection{Relative Frobenius models}

A perfect multiplication pairing gives an algebraic substitute for the
tubular neighbourhood: its double annihilator forces every dual map
on the kernel to take values in that kernel. This is the relative
version of the finite Poincar\'e pairing argument in
\cite[Lemma~4.4]{ParentI}.

\begin{definition}[Relative Frobenius model]\label{def:relative-frobenius}
A relative Frobenius model over a connected CDGA $B$ is a
split CDGA $B\to S\xrightarrow{\pi}B$ that is finite semifree as a
differential $B$-module, together with a closed $B$-linear map
$\lambda:S\to B$ of degree $-N$ for which the multiplication
pairing $s\mapsto\lambda(s\cdot-)$ induces an isomorphism
$S\cong\operatorname{Hom}_B(S,B)[-N]$ of differential $S$-modules,
with the usual shifted action on Hom.
We call $\lambda$ its \emph{relative Frobenius trace}.
\end{definition}
The finiteness is relative to the base: no boundedness assumption is
imposed on its cohomology. The next proposition is the algebraic
step. Its corollary explains how to use the step for a topological
section.

\begin{proposition}[Evaluation in a relative Frobenius model]
\label{prop:relative-frobenius-evaluation}
Let $A\to S\xrightarrow{\pi}A$ be a relative Frobenius model
in the sense of Definition~\ref{def:relative-frobenius}.
Put $J=\ker\pi$. Then canonical biduality
$J\to\mathbb D_S\mathbb D_SJ$ is an isomorphism, and derived
evaluation $J\otimes_S^{\mathbf L}\mathbb D_SJ\to S$ factors
through $J\hookrightarrow S$.
\end{proposition}
\begin{proof}
We first identify the derived duals, and then show that evaluation
preserves $J$.

\emph{Duality.}
The unit and augmentation split $S=A\oplus J$ as differential
$A$-modules. Thus $J$ is a direct summand of a finite semifree
$A$-module. Its $A$-linear dual has the same property: the dual of
a finite semifree module is semifree with the reverse cell order.
Both modules are therefore homotopically projective over $A$, and the natural
map to the double $A$-linear dual is an isomorphism.

The Frobenius isomorphism and coinduction give, for every
differential $S$-module $M$, an isomorphism of complexes
\begin{equation}\label{eq:relative-frobenius-hom}
 \operatorname{Hom}_S(M,S)
 \cong\operatorname{Hom}_A(M,A)[-N].
\end{equation}
Choose a semifree $S$-resolution $P\to J$. Since $S$ is semifree
over $A$, $P$ is semifree over $A$ as well. The quasi-isomorphism
$P\to J$ is consequently an $A$-linear homotopy equivalence.
It follows from \eqref{eq:relative-frobenius-hom} that
\[
 D_J=\operatorname{Hom}_S(J,S)
       \longrightarrow\operatorname{Hom}_S(P,S)
\]
is a quasi-isomorphism. Thus $D_J$ computes $\mathbb D_SJ$.
As an $A$-module, $D_J$ is the shifted dual of $J$, so the same
argument computes $\mathbb D_SD_J$ by ordinary Hom. Applying
\eqref{eq:relative-frobenius-hom} twice, with the graded symmetry
of the multiplication pairing, identifies the natural biduality map
with the double $A$-linear duality isomorphism of $J$.

\emph{The double annihilator.}
The Frobenius pairing identifies $\operatorname{Ann}(J)$ with
the $A$-linear functionals on $S$ that vanish on $J$. Indeed,
$\lambda(sJ)=0$ implies $\lambda(sjb)=0$ for every $j\in J$
and $b\in S$, and perfection of the pairing gives $sJ=0$.
Since $S=A\oplus J$ as $A$-modules, these functionals are exactly
$\operatorname{Hom}_A(S/J,A)\cong A$, generated by $\pi$. Hence $\operatorname{Ann}(J)$ is generated over $A$
by its Frobenius-dual element $\Delta$. Since $J\Delta=0$,
\[
 s\Delta=\pi(s)\Delta\qquad(s\in S),
\]
where coefficients from $A$ act through its unit map. The map
$A\to A\Delta$ is injective, so
\[
 \operatorname{Ann}(\operatorname{Ann}(J))=J.
\]
\emph{Evaluation.}
For any homogeneous $S$-linear map $u:J\to S$,
\[
 a u(x)=(-1)^{|a||u|}u(ax)=0
 \qquad(a\in\operatorname{Ann}(J),\ x\in J).
\]
Thus $u(J)\subset J$. Ordinary evaluation
$J\otimes_S D_J\to S$ therefore takes values in $J$.
Precomposing with semifree replacements of both factors represents
the required factorization of derived evaluation.
\end{proof}

\begin{corollary}[Equality from a relative Frobenius model]
\label{cor:relative-frobenius-equality}
Let $f:X\to Y$ be rationally retractive between simply connected
spaces of finite rational type. Suppose its split minimal model fits into
a strictly commutative square
\[
\begin{tikzcd}[column sep=3.5em,row sep=1.6em]
 R \arrow[r] \arrow[d,"\simeq"'] & A \arrow[d,"\simeq"] \\
 S \arrow[r,"\pi"] & B\rlap{,}
\end{tikzcd}
\]
whose horizontal maps are surjective, and suppose
$B\to S\xrightarrow{\pi}B$ is a relative Frobenius model
(Definition~\ref{def:relative-frobenius}).
Then
\[
 \msecat(f_{\Q})=\secat(f_{\Q})=\relcat(f_{\Q}).
\]
\end{corollary}
\begin{proof}
Proposition~\ref{prop:relative-frobenius-evaluation} gives biduality
and ideal-valued evaluation for $\ker\pi$.
The strictly commutative square gives a quasi-isomorphism
$I\to J=\ker\pi$ by the two short exact sequences of kernels.
The transport lemma \cite[Lemma~4.7]{ParentI}
identifies canonical evaluation and biduality for these two kernels
under extension and restriction along $R\to S$; the underlying derived
equivalence is \cite[Theorem~4.3]{SchwedeShipley2000}.
Its tensor hypothesis holds here: tensoring with a semifree module
preserves quasi-isomorphisms, as follows from its cell filtration
and exactness of tensor products over $\Q$.
Theorem~\ref{thm:imported-evaluation} now gives the conclusion.
\end{proof}

Monic polynomial relations provide a particularly simple source of
relative Frobenius models. The following lemma gives the models for projective bundles and for
the two-sphere factor in the fixed-hyperbolic-fibre family
(Section~\ref{subsec:fixed-hyperbolic-fibre}).

\begin{lemma}[Monic polynomial models]\label{lem:monic-frobenius}
Let $B$ be a connected CDGA and let $x$ have positive even
degree. Suppose that
\[
 p(x)=x^{r+1}+a_r x^r+\cdots+a_1x,\qquad r\geq1,
\]
is homogeneous and that all $a_i\in B$ are closed. Give
$S=B[x]/(p)$ the differential induced by $dx=0$ and the
augmentation $x\mapsto0$. Then $S$ is finite semifree over $B$,
and the coefficient of $x^r$ defines a closed relative Frobenius
trace of degree $-r|x|$.
Moreover, the projection
\[
 (B\otimes\Lambda(x,y),D)\longrightarrow S,
 \qquad Dx=0,\quad Dy=p(x),\quad |y|=(r+1)|x|-1,
\]
which kills $y$, is a quasi-isomorphism compatible with the
augmentations to $B$.
\end{lemma}
\begin{proof}
Monicity gives a differential $B$-module basis
$1,x,\ldots,x^r$. The trace is closed because the coefficients
of $p$ are closed. For its multiplication pairing, the matrix
entry $\lambda(x^{i+j})$ is zero when $i+j<r$ and one when
$i+j=r$. Reversing the columns gives a triangular matrix with
diagonal entries one. This unitriangular matrix is invertible
over the graded base $B$, proving that the pairing is perfect.

Multiplication by $p$ is injective on $B[x]$, even if $B$ has
zero divisors. The kernel of the displayed projection is
$pB[x]\oplus yB[x]$. The degree-$-1$ map
$H(pb)=yb$, $H(yb)=0$, satisfies $DH+HD=\id$ on this kernel.
Thus the projection is a quasi-isomorphism.
\end{proof}

\section{Two families with explicit values}
\label{sec:geometric-examples}

The equality theorems become more informative when the common value
can also be computed. We develop two families. The first consists of
projective bundles over $\mathbb{CP}^k$; already for fibre $S^2$, the
value grows with $k$. The second has a fixed rationally hyperbolic
manifold as fibre and exhibits the same growth. In both cases the
invariant is attached to the section, not to the bundle projection.
We finish with a limitation on smooth realization.

\subsection{Projective bundles and a fixed two-sphere fibre}
\label{sec:projective-examples}

We now compute the family introduced in Section~\ref{sec:introduction}.
Let $k,r\geq1$, let $X=\mathbb{CP}^{k}$, and let $L\to X$ be
the hyperplane line bundle. Write $u=c_1(L)$. A fixed line in the
trivial summand defines a section
\[
 f:X\longrightarrow E=\mathbb P(\underline{\mathbb C}^{\,r}\oplus L)
 \quad\text{of}\quad g:E\longrightarrow X.
\]
Here $\mathbb P$ denotes the bundle of complex lines, so the fibre
of $g$ is $\mathbb{CP}^r$. Let $\mathcal T$ be its tautological
line bundle and put $x=-c_1(\mathcal T)$. Then $g^*u=u$ and
$f^*x=0$. The projective bundle formula
\cite{CMS2010} gives
\begin{equation}\label{eq:projective-cohomology}
 H^*(E;\Q)=
 \Q[u,x]/(u^{k+1},x^{r+1}+ux^r),
 \qquad |u|=|x|=2.
\end{equation}

The cohomology algebra already suggests the answer: $x^{k+r}$ is
nonzero, whereas $x^{k+r+1}$ vanishes. We now prove that this gives
the common value of the three sectional invariants.

\begin{proposition}[The value of the projective section]
\label{prop:projective-section-values}
For the section $f$ just defined,
\[
 \msecat(f_{\Q})=\secat(f_{\Q})=\relcat(f_{\Q})=k+r.
\]
\end{proposition}
\begin{proof}
Theorems~\ref{thm:finite-homotopy-fibre} and~\ref{thm:smooth-section}
already give equality for this family. We use a relative Frobenius
model because the same model also computes the common value from
powers of its augmentation ideal.

\emph{A model that respects the section.}
Set $A=(\Lambda(u,v),d)$, where $|v|=2k+1$, $du=0$ and
$dv=u^{k+1}$. A split minimal model is
\begin{equation}\label{eq:projective-split-model}
 \begin{gathered}
 \varphi:R=A\otimes\Lambda(x,y)\longrightarrow A,
 \qquad \varphi(x)=\varphi(y)=0,\\
 |x|=2,\qquad |y|=2r+1,\qquad
 d_Rx=0,\qquad d_Ry=x^{r+1}+ux^r.
 \end{gathered}
\end{equation}
The inclusion $A\hookrightarrow R$ models $g$.
To verify the compatibility, choose a cocycle $\widehat u$ on $X$
and a primitive for $\widehat u^{\,k+1}$, and pull both back to $E$.
Choose a cocycle $\widehat x$ representing $x$. Since $f^*x=0$,
its restriction to $X$ is exact. More directly, replacing
$\widehat x$ by $\widehat x-g^*f^*\widehat x$ preserves its
cohomology class and makes its restriction strictly zero.
The projective bundle relation makes
$\widehat x^{\,r+1}+g^*\widehat u\,\widehat x^{\,r}$ exact.
If $\eta$ is a primitive, replace it by $\eta-g^*f^*\eta$.
Its differential is unchanged, because the relation restricts to
zero, and this new primitive also restricts to zero. These choices
give a strictly commutative model diagram for $f$ and $g$.

Put $B=\Q[u]/(u^{k+1})$ and
$S=B[x]/(x^{r+1}+ux^r)$, both with zero differential.
Lemma~\ref{lem:monic-frobenius} gives a quasi-isomorphism
\[
 R\longrightarrow A[x]/(x^{r+1}+ux^r).
\]
The subsequent map to $S$ is a quasi-isomorphism as well: as a map
of complexes it is the direct sum of $r+1$ shifts of $A\to B$.
Thus $H(R)$ is exactly the algebra in
\eqref{eq:projective-cohomology}. The map to forms constructed
above sends its generators to $u$ and $x$, hence induces the
required cohomology isomorphism. This verifies that
\eqref{eq:projective-split-model} models the section.

\emph{Equality.}
The projections killing $v$ and $y$ give the commutative square
\[
\begin{tikzcd}[column sep=3.5em,row sep=1.6em]
 R \arrow[r,"\varphi"] \arrow[d,"\simeq"']
   & A \arrow[d,"\simeq"] \\
 S \arrow[r,"\overline\varphi"] & B\rlap{,}
\end{tikzcd}
\]
where $\overline\varphi(u)=u$ and $\overline\varphi(x)=0$.
The coefficient of $x^r$ is a relative Frobenius trace on $S$,
by Lemma~\ref{lem:monic-frobenius}. Consequently
Corollary~\ref{cor:relative-frobenius-equality} identifies the
three invariants.

\emph{The value.}
The $B$-basis $1,x,\ldots,x^r$ gives
\[
 x^{r+j}=(-u)^j x^r\quad(j\geq0),\qquad
 x^{k+r}=(-1)^k u^k x^r\ne0,\qquad x^{k+r+1}=0.
\]
Let $h=k+r$. In $R$, the cocycle $x^h\in I^h$, with $I=(x,y)$,
has nonzero cohomology class. Hence $I^h\hookrightarrow R$ is
nonzero in $\mathcal D(R)$, and the module criterion of
Theorem~\ref{thm:imported-obstruction} gives
$\msecat(f_{\Q})\geq h$.

For the upper bound, put $J=\ker\overline\varphi=(x)$. We have
$J^{h+1}=0$, so the quasi-isomorphism $R\to S$ factors as
\[
 R\longrightarrow R/I^{h+1}\longrightarrow S.
\]
Composing the second map with the inverse of $R\to S$ in the
CDGA homotopy category gives a left inverse of the first. Factor
that quotient as a relative Sullivan inclusion $j:R\hookrightarrow R'$
followed by a quasi-isomorphism. The inclusion also has a left inverse
in the homotopy category.

Lemma~\ref{lem:strict-retraction-representative} gives a strict
retraction of this relative Sullivan factorization.
The quotient criterion for sectional category
\cite[Theorem~3.3]{CV2017} gives $\secat(f_{\Q})\leq h$.
The lower bound and the equality already proved finish the calculation.
\end{proof}

For $r=1$, the fibre is $S^2$ and the common value is $k+1$.
When $k>r$, the base has cohomological dimension $2k>2r$, so this
family lies outside the base-dimension bound of Theorem~\ref{thm:pd-fibre-dimension}.

\subsection{A fixed rationally hyperbolic fibre}\label{subsec:fixed-hyperbolic-fibre}

The next family shows that infinite total rational homotopy of the fibre does not prevent an exact calculation.
Let
\[
 M=(S^3\times S^3)\mathbin{\#}(S^3\times S^3),\qquad
 E_k=\mathbb P(\underline{\mathbb C}\oplus L)
      \longrightarrow\mathbb{CP}^{k},
\]
where $k\geq1$, and let $f_k$ be the section above. Fix $m_0\in M$.
Then
\[
 s_k:\mathbb{CP}^{k}\longrightarrow M\times E_k,
 \qquad b\longmapsto(m_0,f_k(b)),
\]
is a section of the smooth bundle obtained by composing projection
onto $E_k$ with its bundle projection. Its fibre is the fixed
closed eight-dimensional manifold $F=M\times S^2$.

\begin{proposition}[Growing values with a fixed hyperbolic fibre]
\label{prop:hyperbolic-smooth}
The fibre $F$ is simply connected and rationally hyperbolic, and
\[
 \msecat((s_k)_{\Q})=
 \secat((s_k)_{\Q})=\relcat((s_k)_{\Q})=k+3.
\]
For $k\geq5$, the base dimension is greater than the formal
dimension of the fibre.
\end{proposition}
\begin{proof}
The Euler characteristic of $M$ is $-2$. Since a rationally elliptic
space has nonnegative Euler characteristic
\cite[Proposition~32.10]{FHT2001}, $M$ is rationally hyperbolic.
It is a retract of $F$, so $F$ is rationally hyperbolic as well.
The base has cohomological dimension $2k$, whereas $F$ has formal
dimension eight; this gives the last assertion.

It remains to compute the module invariant. The manifold $M$ is
$2$-connected and six dimensional, so the connectivity--dimension
bound gives $\operatorname{cat}(M)\leq2$. A nonzero product of two
classes in $H^3(M;\Q)$ gives the reverse rational lower bound.
Thus $\operatorname{cat}_0(M)=2$, and Hess's theorem
\cite{Hess1991} gives $\msecat(*\to M_{\Q})=2$.
The section $s_k$ is the product of $\{m_0\}\hookrightarrow M$
and $f_k$. Both admit retractions, so module additivity
\cite[Theorem~1]{CPV2018}, applied to fibration replacements of their
rationalizations, and
Proposition~\ref{prop:projective-section-values} give
\[
 \msecat((s_k)_{\Q})
 =\msecat(*\to M_{\Q})+\msecat((f_k)_{\Q})
 =2+(k+1)=k+3.
\]
The smooth-section theorem
(Theorem~\ref{thm:smooth-section}) identifies the other
two invariants with this value.
\end{proof}

For sectional category alone, the product inequality also gives
$\secat((s_k)_\Q)\leq\operatorname{cat}_0(M)+\secat((f_k)_\Q)=k+3$
\cite[Lemma~6.3]{FarberOprea2024}. Together with module additivity,
this recovers the first equality. The smooth-section theorem also
identifies relative category with the same value.

For $k\geq5$, neither Theorem~\ref{thm:finite-homotopy-fibre} nor the stated
base-dimension bound applies to this fibre and base. The family
nevertheless also has a relative Frobenius model. If $C$ is a finite
connected differential Poincar\'e duality model of $M$
\cite{LS2008}, with trace $\lambda_C$ and augmentation $\epsilon_C$,
then over $B=\Q[u]/(u^{k+1})$ it is
\[
 C\otimes B[x]/(x^2+ux)\longrightarrow B,
 \qquad \gamma\otimes P(u,x)\longmapsto\epsilon_C(\gamma)P(u,0).
\]
Splitting the finite complex $C$ into cohomology representatives
and contractible pairs proves finite semifreeness over $B$.
Its trace is $\lambda_C\otimes\operatorname{coeff}_x$;
Lemma~\ref{lem:monic-frobenius} supplies the Frobenius pairing on the
two-sphere factor.
Tensoring the compatible projective model with a model of
$\{m_0\}\hookrightarrow M$ supplies the required compatible square.
Thus the geometric and Frobenius criteria overlap here; the single
calculation in the proposition gives their common value.

\begin{remark}[A Frobenius fibre without a smooth realization]
\label{rem:frobenius-nonsmooth}
The Frobenius and smooth hypotheses differ, even though no example here
requires the Frobenius route alone. Let $C$ be the zero-differential
Poincar\'e algebra with basis $1,a_i,c_i,b_i,\omega$ ($1\leq i\leq4$)
in degrees $0,3,4,5,8$, and nonzero positive-degree products
\[
 a_i b_j=\delta_{ij}\omega,\qquad
 c_i c_j=\delta_{ij}d_i\omega,\qquad
 (d_1,d_2,d_3,d_4)=(1,1,1,2).
\]
Let $F$ realize its minimal Sullivan model. Then $\chi(F)=-2$,
so $F$ is rationally hyperbolic, and its first rational homotopy
occurs in degree three. Its middle pairing has determinant square
class $2$, unchanged by rescaling the orientation since its rank is
four. A closed simply connected eight-manifold has unimodular middle
pairing, of determinant square class $1$ or $-1$; thus none realizes
this rational homotopy type.

Over $X=(\mathbb{CP}^k)_\Q$, $k\geq5$, the constant product section
has relative Frobenius model $B\to B\otimes C\to B$, where
$B=\Q[u]/(u^{k+1})$, with trace $1\otimes\operatorname{coeff}_\omega$.
Its common value is two: the augmentation ideal has cube zero and
$a_1b_1=\omega\ne0$, so the quotient and module criteria give the
two bounds (using Lemma~\ref{lem:strict-retraction-representative}).
Here $N=8$, $c=3$ and $2k\geq10$, outside the base-dimension bound.
Changing the retraction does not change the rational fibre: the shear
$(g,\operatorname{pr}_F):X\times F\to X\times F$ is an isomorphism on
rational homotopy, triangular with identity diagonal blocks.
This separates the hypotheses, but equality also follows from the
constant-product description and Hess's theorem.
\end{remark}

\section{The limit of equality and the remaining questions}
\label{sec:limits}

The four criteria remove the obstruction under different hypotheses.
A strict example shows what can fail and leads to two open questions.

\subsection{What a strict gap detects}

Equality cannot hold for all rationally retractive maps. The
construction in \cite[Theorem~1.1]{Parent2026} gives a retractive
rational map $f:X\to Y$: the source has minimal model
$(\Lambda(a_6,b_{17}),db=a^3)$, and the target is the rationalization
of a simply connected $12$-dimensional CW complex with $11$ cells.
That theorem proves
\[
 \msecat(f)=2<\secat(f)=3.
\]
By $\secat(f)\leq\relcat(f)\leq\msecat(f)+1$
\cite[Theorem~1.3]{ParentI}, one also has $\relcat(f)=3$.
For its split model, Theorem~\ref{thm:imported-obstruction} says
\[
 \iota_*o_3(f)=0,\qquad o_3(f)\ne0.
\]
Choose a semifree resolution $P\to I^3$ and an $R$-linear
nullhomotopy $h:P\to R$ of the inclusion. Its augmentation
$\varphi h:P\to A$ is a degree $-1$ cocycle. The exact sequence
of coefficient complexes gives
\[
 \operatorname{Hom}_{\mathcal D(R)}(I^3,R[-1])
 \longrightarrow
 \operatorname{Hom}_{\mathcal D(R)}(I^3,A[-1])
 \xrightarrow{\delta}
 \operatorname{Hom}_{\mathcal D(R)}(I^3,I),
\]
and $\delta[\varphi h]=o_3(f)$ with the connecting convention
$\delta[z]=[\partial\widetilde z]$: here the lift is $h$, and
$\partial h$ is the inclusion. Thus $[\varphi h]$ is nonzero modulo the image of the first group.
Changing the nullhomotopy changes this class by an element of that
image, so cannot remove it. In particular, the biduality and ideal-valued evaluation
hypotheses of Theorem~\ref{thm:imported-evaluation} cannot both hold
for this kernel. The strict example measures a failure of compatibility
with the ideal. Whether this particular fibre has Poincar\'e duality remains
undetermined here. Deciding it requires the relative minimal model
of a chosen retraction and the cohomology of its fibre algebra
$(\Lambda W,\overline d)$; the two sectional-category values do not
supply those data. In particular, the gap cannot be used to infer
failure of fibre duality without first resolving the question below.

\subsection{Two questions left by the Thom criterion}
\label{sec:open-questions}

\begin{enumerate}[label=\textup{\arabic*.},leftmargin=*]
\item \emph{Does every Poincar\'e duality fibre satisfy derived
self-intersection?} Let $f:X\to Y$ be rationally retractive and let
the fibre of a chosen retraction have finite-dimensional Poincar\'e
duality cohomology. Must its Thom morphism satisfy
$\varphi\tau=\varphi^*(m_e)$ in $\mathcal D(R)$? Equivalently,
is $\kappa_{\mathrm{Th}}$ always zero, or does there exist a section for which
$T_I\ne0$? A positive answer would identify all three sectional
invariants for every such section, without a dimension, finite
homotopy, or smoothness hypothesis. The relative-duality lemma alone
does not answer this question: it proves biduality and leaves precisely
this compatibility undecided.

\item \emph{Can the relative Frobenius criterion give an equality beyond
the other available tests?} Find a rationally retractive map with a
finite relative Frobenius model for which equality is not already
explained by finite total rational homotopy of the fibre, the smooth-section theorem,
the dimension tests, or the constant-product case of Hess's theorem.
Remark~\ref{rem:frobenius-nonsmooth} shows that Frobenius models need
not be manifold models, but its equality also follows from the product
description. The question asks for an example
where the relative pairing supplies genuinely additional information.
\end{enumerate}

\appendix
\section{Relative duality and canonical biduality}
\label{app:relative-duality}

We prove Lemma~\ref{lem:relative-pd-duality}. The finite model is over
$A$, whereas the natural evaluation is $R$-linear. The final step
checks that these two viewpoints identify the same biduality map.

\begin{proof}[Proof of Lemma~\ref{lem:relative-pd-duality}]
\emph{Step 1: a finite module model over the base.}
Choose a contraction $(i_F,p_F,h_F)$ of $(\Lambda W,\overline d)$ onto
$H^*(F;\Q)$, and extend it $A$-linearly. Write
\[
 d_R=d_A\otimes1+1\otimes\overline d+\delta_{\mathrm{tw}}.
\]
Every coefficient of $\delta_{\mathrm{tw}}$ has positive degree in $A$,
hence degree at least two because $A^1=0$. Thus $\delta_{\mathrm{tw}}$ lowers degree in $\Lambda W$
by at least one. The chosen contracting homotopy lowers that degree
by one as well, so the perturbation series are finite on each
element. Choose the usual side conditions
$p_Fi_F=1$, $h_Fi_F=p_Fh_F=h_F^2=0$, and
$dh_F+h_Fd=1-i_Fp_F$. The basic perturbation lemma
\cite[Main Perturbation Lemma~2.4]{Crainic2004} applies. Explicitly, put
$\Omega=\sum_{a\geq0}(-1)^a\delta_{\mathrm{tw}}
(h_F\delta_{\mathrm{tw}})^a$. The transferred data are
$D_{L_F}=d_A+p_F\Omega i_F$,
$i'_F=i_F-h_F\Omega i_F$, $p'_F=p_F-p_F\Omega h_F$, and
$h'_F=h_F-h_F\Omega h_F$.
Our sign convention is $dh_F+h_Fd=1-i_Fp_F$, the negative
of the contracting homotopy convention in the cited perturbation
lemma. With the formulas above, the retract identities and
$\Omega=\delta_{\mathrm{tw}}-\delta_{\mathrm{tw}}h_F\Omega$ give
\[
 p'_Fi'_F=1,\qquad
 d_Rh'_F+h'_Fd_R=1-i'_Fp'_F.
\]
They also give the cochain identities for $i'_F$ and $p'_F$. Thus we obtain an
$A$-module model
\[
 L_F=(A\otimes H^*(F;\Q),D_{L_F})\simeq R.
\]
Its transferred differential, apart from $d_A$, has coefficients
in $A^{\geq2}$ and strictly lowers the degree of a cell in $H^*(F;\Q)$.
Ordering these finitely many cells by degree makes $L_F$ finite
semifree over $A$. Therefore $R$ is $A$-perfect. The cochain
splitting $R=A\oplus I$ makes $I$ a retract of $R$ as an
$A$-module, and hence $A$-perfect too.

\emph{Step 2: fibre integration and its pairing.}
Fix an orientation $H^N(F;\Q)\to\Q$. Its projection on the top-degree
cell of $L_F$, extended as a graded $A$-linear map of degree $-N$,
is a cocycle: the perturbation part of the transferred differential
has no component into that cell. Composing with the transferred projection $R\to L_F$
gives a cochain map $\Pi:R\to A$ of degree $-N$ that restricts to
the orientation on the fibre. The pairing $r\otimes r'\mapsto
\Pi(rr')$ defines
\[
 \Theta_{\Pi}:R\longrightarrow\operatorname{Hom}_A(R,A)[-N],
 \qquad \Theta_{\Pi}(r)(r')=\Pi(rr').
\]
The differential on $\operatorname{Hom}_A(R,A)[-N]$ is
$(-1)^N d_{\operatorname{Hom}}$. Directly,
$d_{\operatorname{Hom}}\Theta_{\Pi}(r)=(-1)^N\Theta_{\Pi}(d_Rr)$;
thus this shifted differential, together with the shifted module
action, makes $\Theta_{\Pi}$ a degree-zero $R$-linear cochain map. To check it by reduction to the
fibre, first transport it to finite semifree $A$-models. Let
$i'_F:L_F\to R$ and $p'_F:R\to L_F$ be the transferred inclusion and
projection. They are inverse $A$-linear homotopy equivalences.
Precomposition with $i'_F$ is therefore a homotopy equivalence
$(i'_F)^*:\operatorname{Hom}_A(R,A)\to\operatorname{Hom}_A(L_F,A)$.
Consequently the cochain map
\[
 \Theta_{\Pi,L_F}=(i'_F)^*[-N]\,\Theta_{\Pi}\,i'_F:
 L_F\longrightarrow\operatorname{Hom}_A(L_F,A)[-N]
\]
represents the underlying $A$-module morphism of $\Theta_{\Pi}$.
The source is finite semifree; the target is too, using the reverse
order on the dual cells. Only on these finite models do we apply
$\Q\otimes_A^{\mathbf L}-$, computed by ordinary tensor product.
Reduction modulo $A^+$ identifies the two sides with
$H^*(F;\Q)$ and $H^*(F;\Q)^\vee[-N]$, and identifies
$\Q\otimes_A\Theta_{\Pi,L_F}$ with the Poincar\'e pairing: the reduced
$i'_F$ chooses cocycle representatives, and the reduced $\Pi$ evaluates
the orientation on their products. This is an isomorphism.

\emph{Step 3: verification of the duality isomorphism.}
Filter both finite models by degree of the coefficient in $A$.
Their finite sets of cell degrees make the filtration finite in
each total degree. The associated graded map is the identity on
$A$ tensored with this Poincar\'e pairing, hence an isomorphism.
Filtered comparison proves that $\Theta_{\Pi,L_F}$, and therefore $\Theta_{\Pi}$,
is a quasi-isomorphism.
This is the relative duality induced by fibre integration; compare
\cite[Proposition~3.2, and Proposition~4.1]{Prigge2026}.

\emph{Step 4: coinduction and biduality.}
A semifree $R$-module resolution is also semifree over $A$, since
$A\to R$ is a relative Sullivan inclusion. Applying derived Hom
to \eqref{eq:relative-pd-isomorphism} and using coinduction gives
\[
 \mathbf R\operatorname{Hom}_R
  \bigl(M,\mathbf R\operatorname{Hom}_A(R,A)\bigr)
 \simeq\mathbf R\operatorname{Hom}_A(M,A),
\]
which proves \eqref{eq:relative-pd-coinduction}.
The coinduction identification sends a homogeneous cochain $h$
to $\Pi\circ h$. Applying it twice gives, for every $R$-module $M$,
\[
 \mathbb D_R\mathbb D_R M
 \simeq
 \mathbf R\operatorname{Hom}_A
    \bigl(\mathbf R\operatorname{Hom}_A(M,A),A\bigr).
\]
The shifts by $N$ and $-N$ cancel. The identification is the one
induced by the pairing $\Theta_{\Pi}$, not an arbitrary isomorphism of the
double duals. It makes the diagram
\[
\begin{tikzcd}[column sep=3em]
 M \arrow[r,"\eta_M^R"] \arrow[dr,"\eta_M^A"'] &
 \mathbb D_R\mathbb D_RM \arrow[d,"\simeq"] \\
 & \mathbf R\operatorname{Hom}_A
 (\mathbf R\operatorname{Hom}_A(M,A),A)
\end{tikzcd}
\]
commute after restriction to $A$. To check the natural map,
choose a semifree resolution $P\to M$. Coinduction and $\Theta_{\Pi}$ give
a quasi-isomorphism
\[
 \Theta_{\Pi,*}:\operatorname{Hom}_R(P,R)
   \longrightarrow\operatorname{Hom}_A(P,A)[-N],
 \qquad h\longmapsto\Pi\circ h.
\]
Choose an $R$-semifree resolution
$q:P'\to\operatorname{Hom}_R(P,R)$. It is also $A$-semifree.
Consequently $\operatorname{Hom}_R(P',R)$ computes the iterated
$R$-derived Hom, while $\Theta_{\Pi,*}q$ identifies $P'[N]$ as an
$A$-semifree resolution of $\operatorname{Hom}_A(P,A)$.
This fixes actual resolutions on both sides of the double-dual
comparison; the first Hom complex itself need not be semifree.

For homogeneous $p\in P$ and $h'\in P'$, the representative of
$\eta_M^R$ on these resolutions sends $p$ to the functional
\[
 h'\longmapsto(-1)^{|p||h'|}\,q(h')(p).
\]
After coinduction, its value is composed with $\Pi$. The corresponding
cochain $k=\Pi q(h')$ in $\operatorname{Hom}_A(P,A)$ has degree
$|k|=|h'|-N$. The canonical identification
$\operatorname{Hom}_A(P',A)[-N]\cong
\operatorname{Hom}_A(P'[N],A)$ contributes $(-1)^{N|p|}$
on a degree-$|p|$ functional. Hence the resulting value is
\[
 (-1)^{|p||h'|+N|p|}\Pi(q(h')(p))
   =(-1)^{|p|(|h'|-N)}k(p)
   =(-1)^{|p||k|}k(p).
\]
This is exactly canonical $A$-linear biduality precomposed with the
resolution $P'[N]\to\operatorname{Hom}_A(P,A)$.
Thus the comparison carries $\eta_M^R$ to $\eta_M^A$; the
cancellation of shifts respects evaluation on the chosen resolutions.
Since $I$ is $A$-perfect, $\eta_I^A$, and hence $\eta_I^R$, is an
isomorphism. Applying \eqref{eq:relative-pd-coinduction} to $M=A$
gives the last assertion.
\end{proof}

\section{The two constructions behind the Thom cocycle}
\label{app:filtered-thom}

The proof of Lemma~\ref{lem:filtered-thom} uses two inputs.
First, its acyclic closure must preserve both positive coefficient
degree and positive relative length. Second, adjoining one relative
generator must induce a quasi-isomorphism on the full Hom complexes.
We prove these two assertions in that order.

The cylinder and the one-generator comparisons are classical.
The point needed here is their simultaneous compatibility with the
relative augmentation and the two filtrations.

\subsection{The relative acyclic closure}

\begin{lemma}\label{lem:relative-thom-cylinder}
Let $A\hookrightarrow R=A\otimes\Lambda W\xrightarrow{\varphi}A$
be the split minimal model \eqref{eq:split-model},
with a finite Sullivan basis $w_1,\ldots,w_k$ ordered by
nondecreasing degree. Minimality and $A^1=0$ make this a Sullivan
order, since every relative generator in $d_Rw_i$ has smaller degree.
Put
$R_{<i}=A\otimes\Lambda(w_1,\ldots,w_{i-1})$.
There is a relative acyclic closure
$\varepsilon:\mathcal C_A=R\otimes\Lambda(sW)\to A$ satisfying
\[
 D(sw_i)=w_i-z_i,\qquad Dz_i=dw_i,\qquad
 z_i\in R_{<i}^+\otimes\Lambda(sw_1,\ldots,sw_{i-1}),
\]
where every monomial of $z_i$ has total relative length at least one
in $W\oplus sW$.  In particular, $D$ does not decrease this length.
\end{lemma}

\begin{proof}
\emph{Step 1: construct the acyclic closure.}
Begin with the algebra
$R^{(1)}\otimes\Lambda(sW\oplus\dot W)$ over $A$, where
$|sw_i|=|w_i|-1$, $|\dot w_i|=|w_i|$, and
\[
 D|_{R^{(1)}}=d_R,\qquad D(sw_i)=\dot w_i,\qquad D(\dot w_i)=0.
\]
Its differential has square zero. Define the degree-$-1$ derivation
$s$ by $s(w_i^{(1)})=sw_i$, $s|_A=s(sW)=s(\dot W)=0$, and put
$\mathcal L=sD+Ds$. Then $s^2=0$, $D\mathcal L=\mathcal L D$, and $s\mathcal L=\mathcal L s$.
The derivation $\mathcal L$ is locally nilpotent: its value on $w_i^{(1)}$
is $\dot w_i$ plus terms involving earlier generators, and it
vanishes on $sW$ and $\dot W$. Induction in the Sullivan order proves
nilpotence on each generator and then on each monomial. Hence
$\exp(\mathcal L)$ is a finite-on-elements CDGA automorphism.

Set $w_i^{(2)}=\exp(\mathcal L)w_i^{(1)}$.
The exponential-cylinder construction
\cite[Proposition~A.4 and Theorem~A.8]{Wakatsuki2016} gives
\[
 D(sw_i)=w_i^{(2)}-w_i^{(1)}
       -\sum_{j\geq1}\frac{(sD)^j(w_i^{(1)})}{j!}.
\]
The change from $\dot w_i$ to $w_i^{(2)}$ is triangular, so it is
invertible and preserves $D^2=0$. The augmentation identifying the
two copies and killing $sW$ is a quasi-isomorphism: in the original
coordinates it removes the contractible pairs $(sw_i,\dot w_i)$.
It is semifree over either copy of $R$.

Specialize the first copy of $W$ to zero. The formula becomes
$D(sw_i)=w_i-z_i$, and $D^2=0$ gives $Dz_i=dw_i$.
Before specialization the cylinder is semifree over the first copy
of $R$, and its augmentation has target free of rank one.
It is therefore a module homotopy equivalence. Tensoring over that
copy with $R\to A$ preserves the homotopy equivalence and yields
$\varepsilon:\mathcal C_A\to A$, with $\varepsilon|_R=\varphi$ and
$\varepsilon(sW)=0$. This is the required acyclic closure.
Every step uses only the current generator and its predecessors.
For $R_i=A\otimes\Lambda(w_1,\ldots,w_i)$, applying the same
construction to this initial segment produces
$\mathcal C_i=R_i\otimes\Lambda(sw_1,\ldots,sw_i)\to A$ with exactly the
same $z_j$ for $j\leq i$. These prefix closures are the compatible
resolutions used in the one-generator induction.

\emph{Step 2: positive coefficient degree.}
The coefficient condition is
$z_i\in R_{<i}^+\otimes\Lambda(sw_1,\ldots,sw_{i-1})$. Before
specialization, its analogue uses the ideal $\mathcal J$ generated
by $A^+$ and the two unsuspended copies of $W$; an element of
$A^+$ counts as one coefficient factor. Relative minimality makes $sD(w_i^{(1)})$
lie in $\mathcal J$: suspension removes at most one factor from a
decomposable relative differential, and vanishes on $A$.
Each further application of $sD$ preserves this property. If $D$
acts on an unsuspended relative factor, relative minimality replaces
it by at least two unsuspended factors before $s$ removes at most one.
If $D$ acts on a suspended factor, it introduces an unsuspended
factor while the original one is retained before $s$ acts.
If the coefficient lies in $A^+$, suspension vanishes on it and on
its differential. Every resulting term therefore still has a
positive coefficient degree in this algebra. Thus specialization gives
$z_i\in R_{<i}^+\Lambda(sw_1,\ldots,sw_{i-1})$.

\emph{Step 3: positive relative length.}
The other condition follows from the augmentation induced by the
cylinder construction: $\varepsilon$ is a cochain map and
$0=\varepsilon D(sw_i)=-\varepsilon(z_i)$. Since
$\ker\varepsilon=(W\oplus sW)$, every monomial of $z_i$ has
positive relative length.

Finally $dW\subset(W)$ and the relative length of each $z_i$ is
positive, so the differential cannot decrease total relative word
length.
\end{proof}

\subsection{The complete filtered comparison}

\begin{lemma}\label{lem:complete-thom-comparison}
Suppose one stage of the preceding construction has the form
\[
 {R_{<i}}=A\otimes\Lambda(w_1,\ldots,w_{i-1}),\quad
 \mathcal C_{<i}={R_{<i}}\otimes\Lambda(sw_1,\ldots,sw_{i-1}),\quad
 R_i={R_{<i}}\otimes\Lambda x,
\]
with $x=w_i$, $t=sx$, and
$\mathcal C_i=R_i\otimes_{R_{<i}}\mathcal C_{<i}\otimes\Lambda t$, where
$Dt=x-z$, $z\in {R_{<i}}^+\mathcal C_{<i}$, and $Dz=dx$.
Use the operator $\Phi$ of
\eqref{eq:thom-phi-odd}--\eqref{eq:thom-phi-even}. The induced map
\[
 \Phi:\operatorname{Hom}_{R_{<i}}(\mathcal C_{<i},{R_{<i}})
       \longrightarrow\operatorname{Hom}_{R_i}(\mathcal C_i,R_i)
\]
is graded ${R_{<i}}$-linear, has degree $\sigma_x=|x|$ for odd $x$ and
$\sigma_x=1-|x|$ for even $x$, and satisfies
$\partial\Phi=-\Phi\partial$.  It induces an isomorphism on
cohomology after shifting by $\sigma_x$.
\end{lemma}

\begin{proof}
\emph{Step 1: the comparison is a map of shifted complexes.}
Both possible values of $\sigma_x$ are odd.  Substitution into the
Hom differential
$\partial \xi=D\xi -(-1)^{|\xi|}\xi D$ verifies
$\partial\Phi=-\Phi\partial$ and graded ${R_{<i}}$-linearity.
For odd $x$, on inputs without $t$ write
$\Phi(\xi)=L_x\xi -(-1)^p\xi L_z$, where $L_v$ means left multiplication
by $v$. The terms involving the old differential cancel because
\[
 L_{dx}\xi=\xi L_{Dz}:
 \qquad \xi((Dz)\omega)=(-1)^{p|dx|}(dx)\xi(\omega)=(dx)\xi(\omega),
\]
since $Dz=dx\in {R_{<i}}$ and $|dx|$ is even. The remaining terms give
$\partial\Phi(\xi)=-\Phi(\partial \xi)$. On positive powers of $t$,
the calculation also uses $Dt=x-z$ and $x^2=z^2=0$.
For even $x$, $t$ is odd, and substituting
$\Phi(\xi)(\omega t)=(-1)^{p+|\omega|}\xi(\omega)$ gives the same identity;
each term $b v$ of $D\omega $ satisfies $|b|+|v|=|\omega|+1$, which fixes
the sign. Thus the comparison respects the full differential,
including that of the base, after the stated shift.

\emph{Step 2: the complete coefficient filtration.}
Filter both Hom complexes by the degree of their ${R_{<i}}$-coefficients.
More precisely, a cochain belongs to $F^a$ when each of its values on
the chosen free module generators has all coefficients in
${R_{<i}}^{\geq a}$; powers of the new generator $x$ in the target are not
counted.  The filtrations start at $F^0$, are decreasing and
separated, and are preserved by the differentials and by $\Phi$.
They are complete coefficient by coefficient.  For a fixed
homogeneous input monomial and a fixed cochain degree, the output
degree is fixed.  More explicitly, in unshifted notation a
${R_{<i}}$-coefficient $b$ in $\xi(v)$ satisfies
\[
 \deg b\leq |v|+|\xi|.
\]
The algebras are nonnegatively graded and $|x|>0$, so only finitely
many degrees of ${R_{<i}}$ can occur in that output.  The same bound applies
to each shifted summand of the cone with its corresponding cochain
degree. This bound is taken separately for each input monomial.

\emph{Step 3: the associated graded comparison.}
In the associated graded complexes, the differential of ${R_{<i}}$, the
differentials of the old module generators, and $dx$ all raise the
filtration and disappear.  So does $z\in {R_{<i}}^+\mathcal C_{<i}$.  The only
remaining module differential is $Dt=x$.  The associated graded map
is therefore the one-generator Koszul comparison, independently at
each old input monomial and each ${R_{<i}}$-coefficient degree.
For an odd exterior generator $x$, the suspended generator $t$ is
polynomial and $D(t^j)=jxt^{j-1}$.  The cohomology class is represented
by the cochain taking $1$ to $x$ and positive powers of $t$ to zero;
the remaining terms contract using the nonzero integers $j$.
For an even polynomial generator $x$, $t$ is exterior and the
comparison extracts the coefficient of $t$, with the sign in
\eqref{eq:thom-phi-even}; the cohomology is the quotient by $(x)$.
Thus both associated graded comparisons induce isomorphisms.
This calculation is entrywise and applies to the full Hom products.

\emph{Step 4: convergence in the full Hom complexes.}
For convergence, consider the mapping cone of the shifted
comparison. It has the same complete filtration, and every
associated graded piece is acyclic.  Given a cocycle $c$, solve for
a primitive modulo $F^1$, subtract its differential, and solve the
remainder modulo $F^2$.  Continuing gives cochains $h_a\in F^a$ such
that
\[
 c-\partial\Bigl(\sum_{a=0}^{n}h_a\Bigr)\in F^{n+1}
 \qquad(n\geq0).
\]
The sum $h=\sum_{a\geq0}h_a$ converges coefficient by coefficient.
Its differential commutes with this limit: for each free module
input, the module differential is a finite linear combination of
inputs, and the comparison term of the cone differential also uses
only finitely many entries, since $z$ is a finite expression.
Each resulting output coefficient stabilizes, so applying the
differential before or after the limit gives the same coefficient.
Separatedness now yields $\partial h=c$.  The cone is acyclic, which
proves that $\Phi$ is a quasi-isomorphism after the shift.
\end{proof}

For each initial segment, the prefix closure from
Lemma~\ref{lem:relative-thom-cylinder} uses the same $z_j$ and resolves
$A$ over $R_i$. Starting with $\operatorname{Hom}_A(A,A)=A$ and
iterating Lemma~\ref{lem:complete-thom-comparison} therefore proves the generator assertion
in Lemma~\ref{lem:filtered-thom}.  The sum of the successive shifts is
\[
 \sum_{w\in\mathcal W_{\mathrm{odd}}}|w|
 -\sum_{w\in\mathcal W_{\mathrm{even}}}(|w|-1)=N.
\]
The independent relative-length control from
Lemma~\ref{lem:relative-thom-cylinder}, combined with the induction in
Lemma~\ref{lem:thom-length}, supplies its ideal-power estimate.

\smallskip
\noindent\emph{Note on citation numbering.}
Numbered citations to \cite{DEH2013} follow the linked arXiv version.
Numbered citations to \cite{ParentI} refer to arXiv:2610.04153v2.

\section*{Competing interests}
The author declares none.

\bigskip
\noindent\textbf{Declaration on the use of generative artificial intelligence.}
During the preparation of this manuscript, the author used ChatGPT
(OpenAI) for language editing, reorganization of the exposition,
\LaTeX{} preparation, bibliographic assistance, and the exploration,
and checking of mathematical arguments. The author also
used Claude (Anthropic) to produce mathematical notes and critical
readings of drafts. The author assumes full responsibility for the
mathematical statements, proofs, references, and conclusions of the article.
\end{document}